\documentclass[article,onecolumn,12pt]{IEEEtran}

\usepackage{cite}
\usepackage{color}
\usepackage{graphicx}
\graphicspath{{fig/}{jpeg/}}
\usepackage[cmex10]{amsmath}
\usepackage{amssymb}
\usepackage{algorithm}
\usepackage{algorithmic}
\usepackage{amsmath,amssymb,lipsum}
\usepackage{hyperref}
\usepackage{comment}
\usepackage{graphicx}
\usepackage[font=small]{caption}
\usepackage{subcaption}
\usepackage{setspace}

\input{mysymbol.sty}
\usepackage{needspace}

\newcommand{\myparagraph}[1]{\needspace{1\baselineskip}\medskip\noindent {\it #1.}}

\renewcommand{\QED}{\hfill\ensuremath{\blacksquare}}
\usepackage{theorem}
\newtheorem{thm}{Theorem}
\newtheorem{lemma}{Lemma}
\newtheorem{proposition}{Proposition}
\newtheorem{corollary}{Corollary}
\newtheorem{assumption}{Assumption}
\newtheorem{remark}{Remark}

\def\Fro{{\mathrm{F}}}
\def\df{{\nabla f}}
\def\tlmin{{\hat{\lambda}_{\min}}}

\def\eps{\epsilon}
\def\bbG{\mathcal{G}}
\date{\today}

\newcommand{\R}{{\mathbb{R}}}

\title{A Decentralized Second-Order Method with Exact Linear Convergence Rate for Consensus Optimization}

\author{Aryan Mokhtari, Wei Shi, Qing Ling, and Alejandro Ribeiro
\thanks{{Work supported by NSF CAREER CCF-0952867, ONR N00014-12-1-0997, and NSFC 61004137. A. Mokhtari and A. Ribeiro are with the Dept. of Electrical and Systems Engineering, University of Pennsylvania, 200 S 33rd St., Philadelphia, PA 19104. Email: \{aryanm, aribeiro\}@seas.upenn.edu. W. Shi is with the Coordinated Science Laboratory, University of Illinois at Urbana-Champaign, 1308 W Main St, Urbana, IL 61801. Email: wilburs@illinois.edu. Q. Ling is with the Dept. of Automation, University of Science and Technology of China, 96 Jinzhao Rd., Hefei, Anhui, 230026, China. Email: qingling@mail.ustc.edu.cn.}}
}

\begin{document}
\thispagestyle{empty}
\maketitle


\begin{abstract}
This paper considers decentralized consensus optimization problems where different summands of a global objective function are available at nodes of a network that can communicate with neighbors only. The proximal method of multipliers is considered as a powerful tool that relies on proximal primal descent and dual ascent updates on a suitably defined augmented Lagrangian. The structure of the augmented Lagrangian makes this problem non-decomposable, which precludes distributed implementations. This problem is regularly addressed by the use of the alternating direction method of multipliers. The exact second order method (ESOM) is introduced here as an alternative that relies on: (i) The use of a separable quadratic approximation of the augmented Lagrangian. (ii) A truncated Taylor's series to estimate the solution of the first order condition imposed on the minimization of the quadratic approximation of the augmented Lagrangian. The sequences of primal and dual variables generated by ESOM are shown to converge linearly to their optimal arguments when the aggregate cost function is strongly convex and its gradients are Lipschitz continuous. Numerical results demonstrate advantages of ESOM relative to decentralized alternatives in solving least squares and logistic regression problems.
\end{abstract}


%
\section{Introduction}

In decentralized consensus optimization problems, components of a global objective function that is to be minimized are available at different nodes of a network. Formally, consider a decision variable $\tbx\in \reals^p$ and a connected network containing $n$ nodes where each node $i$ has access to a local objective function $f_i: \reals^p\to\reals $. Nodes can exchange information with neighbors only and try to minimize the global cost function $\sum_{i=1}^n f_{i}(\tbx)$,
\begin{equation}\label{gen_prob}
\tbx^*\ :=\ \argmin_{\tbx\in\reals^p}\ \sum_{i=1}^n f_{i}(\tbx).
\end{equation}
We assume that the local objective functions $f_i(\tbx)$ are strongly convex. The global objective function $\sum_{i=1}^n f_{i}(\tbx)$, which is the sum of a set of strongly convex functions, is also strongly convex.
Problems like \eqref{gen_prob} arise in decentralized control
\cite{Bullo2009,Cao2013-TII,LopesEtal8}, wireless communication
\cite{Ribeiro10,Ribeiro12}, sensor networks
\cite{Schizas2008-1,KhanEtal10,cRabbatNowak04}, and large scale
machine learning
\cite{bekkerman2011scaling,Tsianos2012-allerton-consensus,Cevher2014}.

Decentralized methods for solving \eqref{gen_prob} can be divided into two classes: primal domain methods and dual domain methods \eqref{gen_prob}. Decentralized gradient descent (DGD) is a well-established primal method that implements gradient descent on a penalized version of \eqref{gen_prob} whose gradient can be separated into per-node components. Network Newton (NN) is a more recent alternative that accelerates convergence of DGD by incorporating second order information of the penalized objective \cite{NN-part1,NN-part2}. Both, DGD and NN, converge to a neighborhood of the optimal argument $\tbx^*$ when using a constant stepsize and converge sublinearly to the exact optimal argument if using a diminishing stepsize. 

Dual domain methods build on the fact that the dual function of \eqref{gen_prob} has a gradient with separable structure. The use of plain dual gradient descent is possible but generally slow to converge \cite{bertsekas1989parallel, ruszczynski2006nonlinear, rabbat2005generalized}. In centralized optimization, better convergence speeds are attained by the method of multipliers (MM) that adds a quadratic augmentation term to the Lagrangian \cite{hestenes1969multiplier, bertsekas2014constrained}, or the proximal (P)MM that adds an additional term to keep iterates close. In either case, the quadratic term that is added to construct the augmented Lagrangian makes distributed computation of primal gradients impossible. This issue is most often overcome with the use of decentralized (D) versions of the alternating direction method of multipliers (ADMM)\cite{Schizas2008-1, BoydEtalADMM11, Shi2014-ADMM}. Besides the ADMM, other methods that use different alternatives to approximate the gradients of the dual function have also been proposed \cite{watanabe1978decomposition, stephanopoulos1975use, mulvey1992diagonal, ruszczynski1995convergence, tappenden2014separable, chatzipanagiotis2013augmented, jakovetic2011cooperative}. The convergence rates of these methods have not been studied except for the DADMM and its variants that are known to converge linearly to the optimal argument when the local functions are strongly convex and their gradients are Lipschitz continuous \cite{Shi2014-ADMM, ling2014dlm, mokhtari2015dqm}. An important observation here is that while all of these methods try to approximate the MM or the PMM, the performance penalty entailed by the approximation has not been studied.

This paper introduces the exact second order method (ESOM) which uses quadratic approximations of the augmented Lagrangians of \eqref{gen_prob} and leads to a set of separable subproblems. Similar to other second order methods, implementation of ESOM requires computation of Hessian inverses. Distributed implementation of this operation is infeasible because while the Hessian of the proximal augmented Lagrangian is neighbor sparse, its inverse is not. ESOM resolves this issue by using the Hessian inverse approximation technique introduced in \cite{ZarghamEtal14, NN-part1, NN-part2}. This technique consists of truncating the Taylor's series of the Hessian inverse to order $K$ to obtain the family of methods ESOM-$K$. Implementation of this expansion in terms of local operations is possible. A remarkable property of all ESOM-$K$ methods is that they can be shown to pay a performance penalty relative to (centralized) PMM that vanishes with increasing iterations.

We begin the paper by reformulating \eqref{gen_prob} in a form more suitable for decentralized implementation (Proposition \ref{new_prop}) and proceed to describe the PMM (Section \ref{sec:problem}). ESOM is a variation of PMM that substitutes the proximal augmented Lagrangian with its quadratic approximation (Section \ref{sec:ESOM}). Implementation of ESOM requires computing the inverse of the Hessian of the proximal augmented Lagrangian. Since this inversion cannot be computed using local and neighboring information, ESOM-$K$ approximates the Hessian inverse with the $K$-order truncation of the Taylor's series expansion of the Hessian inverse. This expansion can be carried out using an inner loop of local operations. This and other details required for decentralized implementation of ESOM-$K$ are discussed in Section \ref{sec:hessian_approx} along with a discussion of how ESOM can be interpreted as a saddle point generalization of the Network Newton methods proposed in \cite{NN-part1,NN-part2} (Remark 1) or a second order version of the EXTRA method proposed in \cite{Shi2014} (Remark 2). 

Convergence analyses of PMM and ESOM are then presented (Section \ref{sec:convg}). Linear convergence of PMM is established (Section \ref{sec:convg_PMM}) and linear convergence factors explicitly derived to use as benchmarks (Theorem \ref{thm:pmm_linear_convg}). In the ESOM analysis (Section \ref{sec:convg_esom}) we provide an upper bound for the error of the proximal augmented Lagrangian approximation (Lemma \ref{esom_approx_error}). We leverage this result to prove linear convergence of ESOM (Theorem \ref{thm:esom_linear_convg}) and to show that ESOM's linear convergence factor approaches the corresponding PMM factor as time grows (Section \ref{sec:rate_comparison}). This indicates that the convergence paths of (distributed) ESOM-$K$ and (centralized) PMM are very close. We also study the dependency of the convergence constant with the algorithm's order $K$. 

ESOM tradeoffs and comparisons with other decentralized methods for solving consensus optimization problems are illustrated in numerical experiments (Section \ref{sec:numerical}) for a decentralized least squares problem (Section \ref{sec:DLS}) and a decentralized logistic regression classification problem (Section \ref{sec:DLR}). Numerical results in both settings verify that larger $K$ leads to faster convergence in terms of number of iterations. However, we observe that all version of ESOM-$K$ exhibit similar convergence rates in terms of the number of communication exchanges. This implies that ESOM-$0$ is preferable with respect to the latter metric and that larger $K$ is justified when computational cost is of interest. Faster convergence relative to EXTRA, Network Newton, and DADMM is observed. We close the paper with concluding remarks (Section \ref{sec:conclusions}).

\myparagraph{\bf Notation} Vectors are written as
$\bbx\in\reals^n$ and matrices as $\bbA\in\reals^{n\times n}$.
Given $n$ vectors $\bbx_i$, the vector
$\bbx=[\bbx_1;\ldots;\bbx_n]$ represents a stacking of the
elements of each individual $\bbx_i$. We use $\|\bbx\|$ and  $\|\bbA\|$ to denote
the Euclidean norm of vector $\bbx$ and matrix $\bbA$, respectively. The norm of vector $\bbx$ with respect to positive definite matrix $\bbA$ is $\|\bbx\|_\bbA:=(\bbx^T\bbA\bbx)^{1/2}$. Given a function $f$ its gradient $\bbx$ is denoted as $\nabla f(\bbx)$ and its Hessian as $\nabla^2 f(\bbx)$. 

%
\section{Proximal method of multipliers}\label{sec:problem}
Let $\bbx_i\in\reals^p$ be a copy of the decision variable $\bbx$ kept at node $i$ and define $\ccalN_{i}$ as the neighborhood of node $i$. Assuming the network is bidirectionally connected, the optimization problem in \eqref{gen_prob} is equivalent to the program 
\begin{align}\label{original_optimization_problem2}
   \{\bbx_i^*\}_{i=1}^n\ := \
   &\argmin_{\{\bbx_{i}\}_{i=1}^n} \ \sum_{i=1}^{n}\ f_{i}(\bbx_{i}), \nonumber\\ 
   &\text{\ s.t.}  \ \bbx_{i}=\bbx_{j}, 
                   \quad \text{for all\ } i, j\in\ccalN_i .
\end{align} 
Indeed, the constraint in \eqref{original_optimization_problem2} enforces the consensus condition $\bbx_1=\dots=\bbx_n$ for any feasible point of \eqref{original_optimization_problem2}. With this condition satisfied, the objective in \eqref{original_optimization_problem2} is equal to the objective function in \eqref{gen_prob} from where it follows that the optimal local variables $\bbx_i^*$ are all equal to the optimal argument $\tbx^*$ of \eqref{gen_prob}, i.e., $\bbx_1^*=\dots=\bbx_n^*=\tbx^*$.

To derive ESOM define $\bbx:=[\bbx_1;\dots;\bbx_n]\in \reals^{np}$ as the concatenation of the local decision variables $\bbx_i$ and the aggregate function $f:\reals^{np}\to\reals$ as $f(\bbx)=f(\bbx_1,\dots,\bbx_n):=\sum_{i=1}^n f_i(\bbx_i)$ as the sum of all the local functions $ f_i(\bbx_i)$. Introduce the matrix $\bbW\in\reals^{n\times n}$ with elements $w_{ij}\geq 0$ representing a weight that node $i$ assigns to variables of node $j$. The weight $w_{ij}=0$ if and only if $j\notin \ccalN_i\cup\{i\}$. The matrix $\bbW$ is further required to satisfy
\begin{equation}\label{eqn_conditions_on_weights}
   \bbW^T=\bbW, \quad 
   \bbW\bbone=\bbone, \quad
   \text{null}(\bbI-\bbW)=\text{span}(\bbone).
\end{equation}
The first condition implies that the weights are symmetric, i.e., $w_{ij}=w_{ji}$. The second condition ensures that the weights of a given node sum up to 1, i.e., $\sum_{j=1}^{n} w_{ij}=1$ for all $i$. Since $\bbW\bbone=\bbone$ we have that $\bbI-\bbW$ is rank deficient. The last condition $\text{null}(\bbI-\bbW)=\text{span}(\bbone)$ makes the rank of $\bbI-\bbW$ exactly equal to $n-1$ \cite{boyd2004fastest}.  

The matrix $\bbW$ can be used to reformulate \eqref{original_optimization_problem2} as we show in the following proposition.

%
\begin{proposition}\label{new_prop}
 Define the matrix $\bbZ:=\bbW\otimes\bbI_p\in\reals^{np} \times \reals^{np}$ as the Kronecker product of the weight matrix $\bbW$ and the identity matrix $\bbI_p$ and consider the definitions of the global vector $\bbx:=[\bbx_1;\dots;\bbx_n]$ and aggregate function $f(\bbx):=\sum_{i=1}^n f_i(\bbx_i)$ . The optimization problem in \eqref{original_optimization_problem2} is equivalent to
\begin{equation}\label{constrained_opt_problem}
   \bbx^* =\argmin_{\bbx \in \reals^{np}} \ f(\bbx)\qquad \text{s.t.} \ (\bbI-\bbZ)^{1/2}\bbx=\bb0.
\end{equation}
I.e., $\bbx^* = [\bbx_1^*;\dots;\bbx_n^*]$ with $\{\bbx_i^*\}_{i=1}^n$ the solution of \eqref{original_optimization_problem2}. 
\end{proposition}
%
%
\begin{myproof} We just show that the constraint $((\bbI_n-\bbW)\otimes\bbI_p)\bbx=(\bbI_{np}-\bbZ)\bbx=\bb0$ is also a consensus constraint. To do so begin by noticing that since $\bbI-\bbW$ is positive semidefinite, $\bbI-\bbZ=(\bbI-\bbW)\otimes\bbI_p$ is also positive semidefinite. Therefore, the null space of the square root matrix $(\bbI-\bbZ)^{1/2}$ is equal to the null space of $\bbI-\bbZ$ and we conclude that satisfying the condition $(\bbI-\bbZ)^{1/2}\bbx$ is equivalent to the consensus condition $\bbx_1=\dots=\bbx_n$. This observation in conjunction with the definition of the aggregate function $f(\bbx)=\sum_{i=1}^n f_i(\bbx_i)$ shows that the programs in \eqref{constrained_opt_problem} and \eqref{eqn_conditions_on_weights} are equivalent. In particular, the optimal solution of \eqref{constrained_opt_problem} is $\bbx^*=[\bbx_1^*;\dots;\bbx_n^*]$ with $\{\bbx_i^*\}_{i=1}^n$ the solution of \eqref{original_optimization_problem2}. \end{myproof}

%
The formulation in \eqref{constrained_opt_problem} is used to define the proximal method of multipliers (PMM) that we consider in this paper. To do so introduce dual variables $\bbv\in\reals^{np}$ to define the augmented Lagrangian $\ccalL(\bbx,\bbv)$ of \eqref{constrained_opt_problem} as
\begin{equation}\label{lagrangian}
\ccalL(\bbx,\bbv)= f(\bbx)+\bbv^T(\bbI-\bbZ)^{1/2}\bbx+\frac{\alpha}{2}\bbx^T(\bbI-\bbZ)\bbx \ \! ,
\end{equation}
where $\alpha$ is a positive constant. Given the properties of the matrix $\bbZ$, the augmentation term $(\alpha/2)\bbx^T(\bbI-\bbZ)\bbx$ is null when the variable $\bbx$ is a feasible solution of \eqref{constrained_opt_problem}. Otherwise, the inner product is positive and behaves as a penalty for the violation of the consensus constraint. 

Introduce now a time index $t\in \mathbb{N}$ and define $\bbx_t$ and $\bbv_t$ as primal and dual iterates at step $t$. The primal variable $\bbx_{t+1}$ is updated by minimizing the sum of the augmented Lagrangian in \eqref{lagrangian} and the proximal term $({\eps}/{2})\|\bbx-\bbx_t\|^2$. We then have that 
\begin{equation}\label{PMM_primal_update}
\bbx_{t+1}=\argmin_{\bbx \in \reals^{np}} \left\{ \ccalL(\bbx,\bbv_t)+\frac{\eps}{2}\|\bbx-\bbx_t\|^2\right\},
\end{equation}
where the proximal coefficient $\eps>0$ is a strictly positive constant. The dual variable $\bbv_t$ is updated by ascending through the gradient of the augmented Lagrangian with respect to the dual variable $\nabla_{\bbv}\ccalL(\bbx_{t+1},\bbv_t)$ with stepsize $\alpha$
\begin{equation}\label{PMM_dual_update}
\bbv_{t+1}=\bbv_{t} +\alpha (\bbI-\bbZ)^{1/2}\bbx_{t+1}.
\end{equation}
The updates in \eqref{PMM_primal_update} and \eqref{PMM_dual_update} for PMM can be considered as a generalization of the method of multipliers (MM), because setting the proximal coefficient $\eps=0$ recovers the updates of MM. The proximal term $({\eps}/{2})\|\bbx-\bbx_t\|^2$ is added to keep the updated variable $\bbx_{t+1}$ close to the previous iterate $\bbx_t$. This does not affect convergence guarantees but improves computational stability. 

The primal update in \eqref{PMM_primal_update} may be computationally costly -- because it requires solving a convex program --  and cannot be implemented in a decentralized manner -- because the augmentation term $({1}/{2\alpha)}\bbx^T(\bbI-\bbZ)\bbx$ in \eqref{lagrangian} is not separable. In the following section we propose an approximation of PMM that makes the minimization in \eqref{PMM_primal_update} computationally economic and separable over nodes of the network. This leads to the set of decentralized updates that define the ESOM algorithm.

%
\section{ESOM: Exact Second-Order Method}\label{sec:ESOM}

To reduce the computational complexity of \eqref{PMM_primal_update} and obtain a separable update we introduce a second order approximation of the augmented Lagrangian in \eqref{lagrangian}. Consider then the second order Taylor's expansion $\ccalL(\bbx,\bbv_t) \approx \ccalL(\bbx_t,\bbv_t) + \nabla_\bbx\ccalL(\bbx_t,\bbv_t)^T(\bbx-\bbx_t) + ({1/2})(\bbx-\bbx_t)^T\nabla^2_\bbx\ccalL(\bbx_t,\bbv_t)(\bbx-\bbx_t)$ of the augmented Lagrangian with respect to $\bbx$ centered around $(\bbx_{t},\bbv_t)$. Using this approximation in lieu of $\ccalL(\bbx,\bbv_t)$ in \eqref{PMM_primal_update} leads to the primal update
\begin{align}\label{ESOM_primal_update}
\bbx_{t+1} = \argmin_{\bbx \in \reals^{np}}  \Big\{\ccalL(\bbx_t,\bbv_t) + \nabla_\bbx\ccalL(\bbx_t,\bbv_t)^T(\bbx-\bbx_t)+ \frac{1}{2}(\bbx-\bbx_t)^T\!\!\left(\nabla^2_\bbx\ccalL(\bbx_t,\bbv_t)+\eps\bbI\right)\!(\bbx-\bbx_t)\Big\}.
\end{align}
The minimization in the right hand side of \eqref{ESOM_primal_update} is of a positive definite quadratic form. Thus, upon defining the Hessian matrix $\bbH_t\in \reals^{np\times np}$ as
\begin{equation}\label{exact_Hessian}
\bbH_t:=\nabla^2f(\bbx_t)+\alpha(\bbI-\bbZ) +\eps\bbI,
\end{equation}
and considering the explicit form of the augmented Lagrangian gradient $\nabla_\bbx\ccalL(\bbx_t,\bbv_t)$ [cf. \eqref{lagrangian}] it follows that the variable $\bbx_{t+1}$ in \eqref{ESOM_primal_update} is given by
\begin{align}\label{ESOM_primal_update_2}
\bbx_{t+1}=\bbx_t - \bbH_t^{-1}\!\!\left[\nabla f(\bbx_t)+ (\bbI-\bbZ)^{1/2}\bbv_t+\alpha(\bbI-\bbZ)\bbx_t\right]\!.
\end{align}
A fundamental observation here is that the matrix $\bbH_t$, which is the Hessian of the objective function in \eqref{ESOM_primal_update}, is block neighbor sparse. By block neighbor sparse we mean that the $(i,j)$th block is non-zero if and only if $j\in\ccalN_i$ or $j=i$. To confirm this claim, observe that $\nabla^2f(\bbx_t)\in \reals^{np\times np}$ is a block diagonal matrix where its $i$th diagonal block is the Hessian of the $i$th local function, $\nabla^2f_i(\bbx_{i,t})\in \reals^{p\times p}$. Additionally, matrix $\eps\bbI_{np}$ is a diagonal matrix which implies that the term $\nabla^2f(\bbx_t)+\eps\bbI_{np}$ is a block diagonal matrix with blocks $\nabla^2f_i(\bbx_{i,t})+\eps\bbI_p$. Further, it follows from the definition of the matrix $\bbZ$ that the matrix $\bbI-\bbZ$ is neighbor sparse. Therefore, the Hessian $\bbH_t$ is also neighbor sparse. Although the Hessian $\bbH_t$ is neighbor sparse, its inverse $\bbH_t^{-1} $ is not. This observation leads to the conclusion that the update in \eqref{ESOM_primal_update_2} is not implementable in a decentralized manner, i.e., nodes cannot implement \eqref{ESOM_primal_update_2} by exchanging information only with their neighbors. 

To resolve this issue, we use a Hessian inverse approximation that is built on truncating the Taylor's series of the Hessian inverse $\bbH_t^{-1} $ as in \cite{NN-part1}. To do so, we try to decompose the Hessian as $\bbH_t=\bbD_t-\bbB$ where $\bbD_t$ is a block diagonal positive definite matrix and $\bbB$ is a neighbor sparse positive semidefinite matrix. In particular, define $\bbD_t$ as 
\begin{equation}
\bbD_t := \nabla^2f(\bbx_t)+\eps\bbI+2\alpha(\bbI-\bbZ_d),
\end{equation}
where $\bbZ_d:=\diag(\bbZ)$. Observing the definitions of the matrices $\bbH_t$ and $\bbD_t$ and considering the relation $\bbB=\bbD_t-\bbH_t$ we conclude that $\bbB$ is given by 
\begin{equation}
\bbB :=\alpha\left(\bbI-2\bbZ_d+\bbZ\right).
\end{equation}
Notice that using the decomposition $\bbH_t=\bbD_t-\bbB$ and by factoring $\bbD_t^{1/2}$, the Hessian inverse can be written as $\bbH_t^{-1}=\bbD_t^{-1/2}(\bbI-\bbD_t^{-1/2}\bbB\bbD_t^{-1/2})^{-1}\bbD_t^{-1/2}$. Observe that the inverse matrix $(\bbI-\bbD_t^{-1/2}\bbB\bbD_t^{-1/2})^{-1}$ can be substituted by its Taylor's series $\sum_{u=0}^{\infty}(\bbD_t^{-1/2}\bbB\bbD_t^{-1/2})^u$; however, computation of the series requires global communication which is not affordable in decentralized settings. Thus, we approximate the Hessian inverse $\bbH_t^{-1}$ by truncating the first $K+1$ terms of its Taylor's series which leads to the Hessian inverse approximation $\tbH_t^{-1}(K)$,
\begin{equation}\label{Hessian_approx}
\tbH_t^{-1}(K):=\bbD_t^{-1/2}\
\sum_{u=0}^K\left(\bbD_t^{-1/2}\bbB\bbD_t^{-1/2}\right)^u\
\bbD_t^{-1/2}.
\end{equation}
Notice that the approximate Hessian inverse $\tbH_t^{-1}(K)$ is $K$-hop block neighbor sparse, i.e., the $(i,j)$th block is nonzero if and only if there is at least one path between nodes $i$ and $j$ with length $K$ or smaller. 

We introduce the Exact Second-Order Method (ESOM) as a second order method for solving decentralized optimization problems which substitutes the Hessian inverse in update \eqref{ESOM_primal_update_2} by its $K$ block neighbor sparse approximation $\hbH_k^{-1}(K)$ defined in \eqref{Hessian_approx}. Therefore, the primal update of ESOM is
\begin{align}\label{ESOM_primal_update_3}
\bbx_{t+1}=\bbx_t - \tbH_t^{-1}(K)\!\!\left[\nabla f(\bbx_t)+ (\bbI-\bbZ)^{1/2}\bbv_t+\alpha(\bbI-\bbZ)\bbx_t\right]\!.
\end{align}
The ESOM dual update is identical to the update in \eqref{PMM_dual_update},
\begin{equation}\label{ESOM_dual_update}
\bbv_{t+1}=\bbv_{t} +\alpha (\bbI-\bbZ)^{1/2}\bbx_{t+1}.
\end{equation}
Notice that ESOM is different from PMM in approximating the augmented Lagrangian in the primal update of PMM by a second order approximation. Further, ESOM approximates the Hessian inverse of the augmented Lagrangian by truncating the Taylor's series of the Hessian inverse which is not necessarily neighbor sparse. In the following subsection we study the implantation details of the updates in \eqref{ESOM_primal_update_3} and \eqref{ESOM_dual_update}.

\subsection{Decentralized implementation of ESOM}\label{sec:hessian_approx}

The updates in \eqref{ESOM_primal_update_3} and \eqref{ESOM_dual_update} show that ESOM is a second order approximation of PMM. Although these updates are necessary for understanding the rationale behind ESOM, they are not implementable in a decentralized fashion since the matrix $(\bbI-\bbZ)^{1/2}$ is not neighbor sparse. To resolve this issue, define the sequence of variables $\bbq_t$ as $\bbq_t:=(\bbI-\bbZ)^{1/2}\bbv_t$. Considering the definition of $\bbq_t$, the primal update in \eqref{ESOM_primal_update_3} can be written as
\begin{align}\label{ESOM_primal_update_4}
\bbx_{t+1}=\bbx_t - \tbH_t^{-1}(K) \big(\nabla f(\bbx_t)+\bbq_t+{\alpha}(\bbI-\bbZ)\bbx_t\big).
\end{align}
By multiplying the dual update in \eqref{ESOM_dual_update} by $(\bbI-\bbZ)^{1/2}$ from the left hand side and using the definition $\bbq_t:=(\bbI-\bbZ)^{1/2}\bbv_t$ we obtain that 
\begin{equation}\label{ESOM_dual_update_4}
\bbq_{t+1}=\bbq_t +\alpha (\bbI-\bbZ)\bbx_{t+1}.
\end{equation}
Notice that the system of updates in \eqref{ESOM_primal_update_4} and \eqref{ESOM_dual_update_4} is equivalent to the updates in  \eqref{ESOM_primal_update_3} and \eqref{ESOM_dual_update}, i.e., the sequences of variables $\bbx_t$ generated by them are identical. Nodes can implement the primal-dual updates in \eqref{ESOM_primal_update_4} and \eqref{ESOM_dual_update_4} in a decentralized manner, since the squared root matrix $(\bbI-\bbZ)^{1/2}$ is eliminated from the updates and nodes can compute the products $(\bbI-\bbZ)\bbx_t$ and $(\bbI-\bbZ)\bbx_{t+1}$ by exchanging information with their neighbors. 

To characterize the local update of each node for implementing the updates in \eqref{ESOM_primal_update_4} and \eqref{ESOM_dual_update_4}, define
\begin{equation}\label{gradient}
\bbg_t:=\nabla_{\bbx} \ccalL(\bbx_t,\bbv_t)= \nabla f(\bbx_t)+\bbq_t+{\alpha}(\bbI-\bbZ)\bbx_t,
\end{equation}
as the gradient of the augmented Lagrangian in \eqref{lagrangian}. Further, define the primal descent direction $\bbd_{t}(K)$ with $K$ levels of approximation as 
\begin{equation}
\bbd_t(K):=-\tbH_t^{-1}(K)\ \! \bbg_t,
\end{equation}
which implies that the update in \eqref{ESOM_primal_update_4} can be written as $\bbx_{t+1}=\bbx_t+\bbd_t(K)$. Based on the mechanism of the Hessian inverse approximation in \eqref{Hessian_approx}, the descent directions $\bbd_{t}(k)$ and $\bbd_{t}(k+1)$ satisfy the condition
\begin{equation}\label{relation}
\bbd_{t}{(k+1)}= \bbD_t^{-1}\bbB \bbd_{t}{(k)}-\bbD_t^{-1}\bbg_t.
\end{equation}
Define $\bbd_{i,t}{(k)}$ as the descent direction of node $i$ at step $t$ which is the $i$th element of the global descent direction $\bbd_{t}{(k)}=[\bbd_{1,t}{(k)};\dots;\bbd_{n,t}{(k)}]$. Therefore, the localized version of the relation in \eqref{relation} at node $i$ is given by
\begin{equation}\label{relation_local}
\bbd_{i,t}{(k+1)}= \bbD_{ii,t}^{-1}\sum_{j=i,j\in\ccalN_i} \bbB_{ij} \bbd_{j,t}{(k)}-\bbD_{ii,t}^{-1}\bbg_{i,t}.
\end{equation}
The update in \eqref{relation_local} shows that node $i$ can compute its $(k+1)$th descent direction $\bbd_{i,t}{(k+1)}$ if it has access to the $k$th descent direction $\bbd_{i,t}(k)$ of itself and its neighbors $\bbd_{j,t}(k)$ for $j\in\ccalN_i$. Thus, if nodes initialize with the ESOM-$0$ descent direction $\bbd_{i,t}{(0)}=-\bbD_{ii,t}^{-1}\bbg_{i,t}$ and exchange their descent directions with their neighbors for $K$ rounds and use the update in \eqref{relation_local}, they can compute their local ESOM-$K$ descent direction $\bbd_{i,t}{(K)}$. Notice that the $i$th diagonal block $\bbD_{t}$ is given by $\bbD_{ii,t}:=\nabla^2 f_{i}(\bbx_{i,t}) + (2\alpha(1-w_{ii})+\eps)\bbI$, where $\bbx_{i,t}$ is the primal variable of node $i$ at step $t$. Thus, the block $\bbD_{ii,t}$ is locally available at node $i$. Moreover, node $i$ can evaluate the blocks $\bbB_{ii}=\alpha(1-w_{ii})\bbI$ and $\bbB_{ij}=\alpha w_{ij}\bbI$ without extra communication. In addition, nodes can compute the gradient $\bbg_t$ by communicating with their neighbors. To confirm this claim observe that the $i$th element of $\bbg_t=[\bbg_{1,t};\dots;\bbg_{n,t}]$ associated with node $i$ is given by
\begin{equation}\label{local_gradient}
\bbg_{i,t}:=\nabla f_i(\bbx_{i,t})+\bbq_{i,t}+{\alpha}(1-w_{ii})\bbx_{i,t}-\alpha\sum_{j\in \ccalN_i}w_{ij}\bbx_{j,t},
\end{equation}
where $\bbq_{i,t}\in \reals^{p}$ is the $i$th element of $\bbq_t=[\bbq_{1,t};\dots;\bbq_{n,t}]$ and $\bbx_{i,t}$ the primal variable of node $i$ at step $t$ and they are both available at node $i$. Hence, the update in \eqref{ESOM_primal_update_4} can be implemented in a decentralized manner. Likewise, nodes can implement the dual update in \eqref{ESOM_dual_update_4} using the local update
\begin{equation}\label{ESOM_dual_update_local}
\bbq_{i,t+1}=\bbq_{i,t} +\alpha (1-w_{ii})\bbx_{i,t+1}- \alpha \sum_{j\in \ccalN_i}w_{ij}\bbx_{j,t+1},
\end{equation}
which requires access to the local primal variable $\bbx_{j,t+1}$ of the neighboring nodes $j\in \ccalN_{i}$.

%
\begin{algorithm}[t]{\small
\caption{ESOM-$K$ method at node $i$}\label{algo_ESOM} 
\begin{algorithmic}[1] {
\REQUIRE  Initial iterates $\bbx_{i,0}={\bbx_{j,0}}=\bb0\ {\forall\ j\in \ccalN_i}$ and $\bbq_{i,0}=\bb0$. 
\STATE $\bbB$ blocks: $\bbB_{ii}=\alpha(1-w_{ii})\bbI$ and $\bbB_{ij}=\alpha w_{ij}\bbI$
\FOR {$t=0,1,2,\ldots$}
   \STATE $\bbD$ block: $\displaystyle{\bbD_{ii,t}= \nabla^2 f_{i}(\bbx_{i,t}) + (2\alpha (1-w_{ii})+\eps)\bbI }$
   \STATE Compute  
          $\displaystyle{
          \bbg_{i,t} \!=\! \nabla f_{i}(\bbx_{i,t})}+{\bbq_{i,t}}+ \alpha{(1-w_{ii})}\bbx_{i,t} 
                       -{\alpha} \sum_{j\in \mathcal{N}_i} {w_{ij}} \bbx_{j,t}$  
   \STATE Compute ESOM-0 descent direction $\bbd_{i,t}{(0)}=-\bbD_{ii,t}^{-1}\bbg_{i,t}$\\ 
   \FOR  {$k=  0, \ldots, K-1$ } 
      \STATE Exchange $\bbd_{i,t}{(k)}$ with neighbors $j\in \ccalN_i$
      \STATE Compute
             $\displaystyle{  
             \bbd_{i,t}{(k+1)}\! = \bbD_{ii,t}^{-1}
             			\bigg[
                                   \sum_{j\in \mathcal{N}_i,j=i}\!\!\! \bbB_{ij} \bbd_{j,t}{(k)} 
                                    - \bbg_{i,t}\bigg]}$        
   \ENDFOR
          %
          \STATE Update primal iterate: 
          $\displaystyle{\bbx_{i,t+1}=\bbx_{i,t} + \bbd_{i,t}{(K)}}$.
             \STATE Exchange iterates $\bbx_{i,t}$ with neighbors $\displaystyle{j\in \mathcal{N}_i}$.
              \STATE Update dual iterate: \\
          $\displaystyle{\bbq_{i,t+1}=\bbq_{i,t} +\alpha (1-w_{ii})\bbx_{i,t+1}- \alpha \sum_{j\in \ccalN_i}w_{ij}\bbx_{j,t+1}}$.
\ENDFOR}
\end{algorithmic}}\end{algorithm}
%

The steps of ESOM-$K$ are summarized in Algorithm \ref{algo_ESOM}. The core steps are Steps 5-9 which correspond to computing the ESOM-$K$ primal descent direction $\bbd_{i,t}(K)$. In Step 5, Each node computes its initial descent direction $\bbd_{i,t}(0)$ using the block $\bbD_{ii,t}$ and the local gradient $\bbg_{i,t}$ computed in Steps 3 and 4, respectively. Steps 7 and 8 correspond to the recursion in \eqref{relation_local}. In step 7, nodes exchange their $k$th level descent direction $\bbd_{i,t}(k)$ with their neighboring nodes to compute the $(k+1)$th descent direction $\bbd_{i,t}(k+1)$ in Step 8. The outcome of this recursion is the $K$th level descent direction $\bbd_{i,t}(K)$ which is required for the update of the primal variable $\bbx_{i,t}$ in Step 10. Notice that the blocks of the neighbor sparse matrix $\bbB$, which are required for step 8, are computed and stored in Step 1. After updating the primal variables in Step 10, nodes exchange their updated variables $\bbx_{i,t+1}$ with their neighbors $j\in \ccalN_i$ in Step 11. By having access to the decision variable of neighboring nodes, nodes update their local dual variable $\bbq_{i,t}$ in Step 12.

\begin{remark} \normalfont
The proposed ESOM algorithm solves problem \eqref{constrained_opt_problem} in the dual domain by defining the proximal augmented Lagrangian. It is also possible to solve problem \eqref{constrained_opt_problem} in the primal domain by solving a penalty version of \eqref{constrained_opt_problem}. In particular, by using the quadratic penalty function $(1/2) \|.\|^2$ for the constraint $(\bbI-\bbZ)^{1/2}\bbx$ with penalty coefficient $\alpha$, we obtain the penalized version of \eqref{constrained_opt_problem}
\begin{equation}\label{penalized_function}
\hbx^*\ :=\ \argmin_{\bbx\in \reals^{np}} f(\bbx)+\frac{\alpha}{2}\bbx^T(\bbI-\bbZ)\bbx,
\end{equation}
where $\hbx^*$ is the optimal argument of the penalized objective function. Notice that $\hbx^*$ is not equal to the optimal argument $\bbx^*$ and the distance $\|\bbx^*-\hbx^*\|$ is in the order of $O(1/\alpha)$. The objective function in \eqref{penalized_function} can be minimized by descending through the gradient descent direction which leads to the update of decentralized gradient descent (DGD) \cite{Nedic2009}. The convergence of DGD can be improved by using Newton's method. Notice that the Hessian of the objective function in \eqref{penalized_function} is given by 
\begin{equation}\label{penalty_Hessian}
\hbH:=\nabla^2f(\bbx)+\alpha(\bbI-\bbZ).
\end{equation}
The Hessian $\hbH$ in \eqref{penalty_Hessian} is identical to the Hessian $\bbH$ in \eqref{exact_Hessian} except for the term $\eps\bbI$. Therefore, the same technique for approximating the Hessian inverse $\hbH^{-1}$ can be used to approximate the Newton direction of the penalized objective function in \eqref{penalized_function} which leads to the update of the Network Newton (NN) methods  \cite{NN-part1,NN-part2}. Thus, ESOM and NN use an approximate decentralized variation of Newton's method for solving two different problems. In other words, ESOM uses the approximate Newton direction for minimizing the augmented Lagrangian of \eqref{constrained_opt_problem}, while NN solves a penalized version of \eqref{constrained_opt_problem} using this approximation. This difference justifies the reason that the sequence of iterates generated by ESOM converges to the optimal argument $\bbx^*$ (Section \ref{sec:convg}), while NN converges to a neighborhood of $\bbx^*$.

\end{remark}

\begin{remark} \normalfont
ESOM approximates the augmented Lagrangian $\ccalL(\bbx,\bbv)$ in \eqref{PMM_primal_update} by its second order approximation. If we substitute the augmented Lagrangian by its first order approximation we can recover the update of EXTRA proposed in \cite{Shi2014}. To be more precise, we can substitute $\ccalL(\bbx,\bbv_t)$ by its first order approximation $\ccalL(\bbx_t,\bbv_t) + \nabla\ccalL(\bbx_t,\bbv_t)^T(\bbx-\bbx_t)$ near the point $(\bbx_t,\bbv_t)$ to update the primal variable $\bbx$. Considering this substitution and the definition of the augmented Lagrangian  
in \eqref{lagrangian} It follows that the update for the primal variable $\bbx$ can be written as
\begin{equation}\label{EXTRA_primal_update_2}
\bbx_{t+1}=\bbx_t - \frac{1}{\eps} \left[\nabla f(\bbx_t)+ {\alpha}(\bbI-\bbZ)^{1/2}\bbv_t+{\alpha}(\bbI-\bbZ)\bbx_t\right].
\end{equation}
By subtracting the update at step $t-1$ from the update at step $t$ and using the dual variables relation that $\bbv_{t+1}=\bbv_{t} +\alpha (\bbI-\bbZ)^{1/2}\bbx_{t+1}$ we obtain the update
\begin{align}\label{EXTRA_primal_update_3}
\bbx_{t+1}
= \left( 2\bbI -\left[\frac{\alpha}{\eps}+\frac{\alpha^2}{\eps}\right](\bbI-\bbZ)\right)\bbx_t- \left( \bbI -\frac{\alpha}{\eps}(\bbI-\bbZ) \right)\bbx_{t-1}
	-\frac{1}{\eps}(\nabla f(\bbx_t)-\nabla f(\bbx_{t-1})).
\end{align}
The update in \eqref{EXTRA_primal_update_3} shows a first-order approximation of the PMM. It is not hard to show that for specific choices of $\alpha$ and $\eps$, the update in \eqref{EXTRA_primal_update_3} is equivalent to the update of EXTRA in \cite{Shi2014}. Thus, we expect to observe faster convergence for ESOM relative to EXTRA as it incorporates second-order information. This advantage is studied in Section \ref{sec:numerical}.

\end{remark}

\section{Convergence Analysis}\label{sec:convg}

In this section we study convergence rates of PMM and ESOM. First, we show that the sequence of iterates $\bbx_t$ generated by PMM converges linearly to the optimal argument $\bbx^*$. Although, PMM cannot be implemented in a decentralized fashion, its convergence rate can be used as a benchmark for evaluating performance of ESOM. We then follow the section by analyzing convergence properties ESOM. We show that ESOM exhibits a linear convergence rate and compare its factor of linear convergence with the linear convergence factor of PMM. In proving these results we consider the following assumptions.

%
\begin{assumption}\label{convexity_assumption} 
The local objective functions $f_i(\bbx)$ are twice differentiable and the eigenvalues of the local objective functions Hessian $\nabla^2 f(\bbx)$ are bounded by positive constants $0<m\leq M<\infty$, i.e. 
\begin{equation}\label{local_hessian_eigenvlaue_bounds00}
m\bbI  \ \preceq\ \nabla^2 f_i(\bbx_i)\ \preceq\ M\bbI, 
\end{equation}
for all $\bbx_i \in \reals^p$ and $i=1,\dots,n$.
\end{assumption}


The lower bound in \eqref{local_hessian_eigenvlaue_bounds00} is equivalent to the condition that the local objective functions $f_i$ are strongly convex with constant $m>0$. The upper bound for the eigenvalues of the Hessians $\nabla^2 f_i$ implies that the gradients of the local objective functions $\nabla f_i$ are Lipschitz continuous with constant $M$. Notice that the global objective function $\nabla^2 f(\bbx)$ is a block diagonal matrix where its $i$th diagonal block is $\nabla^2 f_i(\bbx_i)$. Therefore, the bounds on the eigenvalues of the local Hessians $\nabla^2 f_i(\bbx_i)$ in \eqref{local_hessian_eigenvlaue_bounds00} also hold for the global objective function Hessian $\nabla^2 f(\bbx)$. I.e.,
\begin{equation}\label{local_hessian_eigenvlaue_bounds}
m\bbI \ \preceq \ \nabla^2 f(\bbx) \ \preceq\  M\bbI,
\end{equation}
for all $\bbx \in \reals^{np}$. Thus, the global objective function $f$ is also strongly convex with constant $m$ and its gradients $\nabla f$ are Lipschitz continuous with constant $M$.

\subsection{Convergence of Proximal Method of Multipliers (PMM)}\label{sec:convg_PMM}

Convergence rate of PMM can be considered as a benchmark for the convergence rate of ESOM. To establish linear convergence of PMM, We first study the relationship between the primal $\bbx$ and dual $\bbv$ iterates generated by PMM and the optimal arguments $\bbx^*$ and $\bbv^*$ in the following lemma.

\begin{lemma}\label{lemma_pmm}
Consider the updates for the proximal method of multipliers in \eqref{PMM_primal_update} and \eqref{PMM_dual_update}. The sequences of primal and dual iterates generated by PMM satisfy
\begin{align}\label{opt_res_pmm1}
\bbv_{t+1}-\bbv_{t} -\alpha (\bbI-\bbZ)^{1/2}(\bbx_{t+1}-\bbx^*)=\bb0,
\end{align}
 and
\begin{align}\label{opt_res_pmm2}
\nabla f(\bbx_{t+1})-\nabla f(\bbx^*) &+(\bbI-\bbZ)^{1/2}(\bbv_{t+1}-\bbv^*)
+\eps (\bbx_{t+1}-\bbx_t)=\bb0.
\end{align}
\end{lemma}

\begin{myproof}
See Appendix \ref{app_lemma_pmm}.
\end{myproof}

Considering the preliminary results in \eqref{opt_res_pmm1} and \eqref{opt_res_pmm2}, we can state convergence results of PMM. To do so, we prove linear convergence of a Lyapunov function of the primal $\|\bbx_t-\bbx^*\|^2$ and dual $\|\bbv_t-\bbv^*\|^2$ errors. To be more precise, we define the vector $\bbu\in \reals^{2np}$ and matrix $\bbG\in \reals^{np\times np}$ as
\begin{equation}\label{def:u_G}
\bbu=\left[
       \begin{array}{c}
         \bbv \\
         \bbx \\
       \end{array}
     \right],\
\bbG=\left[
      \begin{array}{cc}
        \bbI & 0 \\
        0 &\alpha \eps\bbI \\
      \end{array}
    \right].
\end{equation}
Notice that the sequence $\bbu_t$ is the concatenation of the dual variable $\bbv_t$ and primal variable $\bbx_t$. Likewise, we can define $\bbu^*$ as the concatenation of the optimal arguments $\bbv^*$ and $\bbx^*$. We proceed to prove that the sequence $\|\bbu_t-\bbu^*\|_\bbG^2$ converges linearly to null. Observe that $\|\bbu_t-\bbu^*\|_\bbG^2$ can be simplified as $\|\bbv_t-\bbv^*\|^2+\alpha\eps\|\bbx_t-\bbx^*\|^2$. This observation shows that $\|\bbu_t-\bbu^*\|_\bbG^2$ is a Lyapunov function of the primal $\|\bbx_t-\bbx^*\|^2$ and dual $\|\bbv_t-\bbv^*\|^2$ errors. Therefore, linear convergence of the sequence $\|\bbu_t-\bbu^*\|_\bbG^2$ implies linear convergence of the sequence $\|\bbx_t-\bbx^*\|^2$. In the following theorem, we show that the sequence $\|\bbu_t-\bbu^*\|_\bbG^2$ converges to zero at a linear rate.

\begin{thm}\label{thm:pmm_linear_convg}
Consider the proximal method of multipliers as introduced in \eqref{PMM_primal_update} and \eqref{PMM_dual_update}. Consider $\beta>1$ as an arbitrary constant strictly larger than $1$ and define $\tlmin(\bbI-\bbZ)$ as the smallest non-zero eigenvalue of the matrix $\bbI-\bbZ$. Further, recall the definitions of the vector $\bbu$ and matrix $\bbG$ in \eqref{def:u_G}. If Assumption \ref{convexity_assumption} holds, then the sequence of Lyapunov functions $\|\bbu_{t}-\bbu^*\|_\bbG^2$ generated by PMM satisfies 
\begin{equation}\label{proof_0}
\|\bbu_{t+1}-\bbu^*\|_\bbG^2 \ \leq\  \frac{1}{1+\delta} \ \|\bbu_{t}-\bbu^*\|_\bbG^2,
\end{equation}
where the constant $\delta$ is given by
\begin{align}\label{pmm_delta}
\delta=\min \Bigg\{ 	
	\frac{2\alpha\tlmin(\bbI-\bbZ)}{\beta(m+M)},\frac{2 m M}{\eps(m+M)} 
	,\frac{(\beta-1)\alpha\tlmin(\bbI-\bbZ)}{\beta\eps} \Bigg\}.
\end{align}
\end{thm}

\begin{myproof}
See Appendix \ref{app_theorem_pmm}.
\end{myproof}

The result in Theorem \ref{thm:pmm_linear_convg} shows linear convergence of the sequence $\|\bbu_{t}-\bbu^*\|_\bbG^2$ generated by PMM where the factor of linear convergence is $1/(1+\delta)$. Observe that larger $\delta$ implies smaller linear convergence factor $1/(1+\delta)$ and faster convergence. Notice that all the terms in the minimization in \eqref{pmm_delta} are positive and therefore the constant $\delta$ is strictly larger than 0. In addition, the result in Theorem \ref{thm:pmm_linear_convg} holds for any feasible set of parameters $\beta>1$, $\eps>0$, and $\alpha>0$; however, maximizing the parameter $\delta$ requires properly choosing the set of parameters $\beta$, $\eps$, and $\alpha$. 

Observe that when the first positive eigenvalue $\tlmin(\bbI-\bbZ)$ of the matrix $\bbI-\bbZ$ , which is the second smallest eigenvalue of $\bbI-\bbZ$, is small the constant $\delta$ becomes close to zero and convergence becomes slow. Notice that small $\tlmin(\bbI-\bbZ)$ shows that the graph is not highly connected. This observation matches the intuition that when the graph has less edges the speed of convergence is slower. Additionally, the upper bounds in \eqref{pmm_delta} show that when the condition number $M/m$ of the global objective function $f$ is large, $\delta$ becomes small and the linear convergence becomes slow. 

Although PMM enjoys a fast linear convergence rate, each iteration of PMM requires infinite rounds of communications which makes it infeasible. In the following section, we study convergence properties of ESOM as a second order approximation of PMM that is implementable in decentralized settings. 

\subsection{Convergence of ESOM}\label{sec:convg_esom}

We proceed to show that the sequence of iterates $\bbx_t$ generated by ESOM converges linearly to the optimal argument $\bbx^*=[\tbx^*;\dots;\tbx^*]$. To do so, we first prove linear convergence of the Lyapunov function $\|\bbu_{t}-\bbu^*\|_\bbG^2$ as defined in \eqref{def:u_G}. Moreover, we show that by increasing the Hessian inverse approximation accuracy, ESOM factor of linear convergence can be arbitrary close to the linear convergence factor of PMM in Theorem \ref{thm:pmm_linear_convg}.

Notice that ESOM is built on a second order approximation of the proximal augmented Lagrangian used in the update of PMM. To guarantee that the second order approximation suggested in ESOM is feasible, the local objective functions $f_i$ are required to be twice differentiable as assumed in Assumption \ref{convexity_assumption}. The twice differentiability of the local objective functions $f_i$ implies that the aggregate function $f$, which is the sum of a set of twice differentiable functions, is also twice differentiable. This observation shows that the global objective function $\nabla^2f(\bbx)$ is definable. Considering this observation, we prove some preliminary results for the iterates generated by ESOM in the following lemma.

\begin{lemma}\label{lemma_esom}
Consider the updates of ESOM in \eqref{ESOM_primal_update_3} and \eqref{ESOM_dual_update}. Recall the definitions of the augmented Lagrangian  Hessian $\bbH_t$ in \eqref{exact_Hessian} and the approximate Hessian inverse $\tbH_t^{-1}(K)$ in \eqref{Hessian_approx}. If Assumption \ref{convexity_assumption} holds, then the primal and dual iterates generated by ESOM satisfy
\begin{align}\label{opt_res_ESOM1}
\bbv_{t+1}-\bbv_{t} -\alpha(\bbI-\bbZ)^{1/2}(\bbx_{t+1}-\bbx^*)=\bb0.
\end{align}
 Moreover, we can show that
\begin{align}\label{opt_res_ESOM2}
&\nabla f(\bbx_{t+1})-\nabla f(\bbx^*) + (\bbI-\bbZ)^{1/2}(\bbv_{t+1}-\bbv^*) +\eps (\bbx_{t+1}-\bbx_t)+\bbe_t=\bb0,
\end{align}
where the error vector $\bbe_t$ is defined as 
\begin{align}\label{esom_error_vec}
\bbe_t &:= \nabla f(\bbx_t)+\nabla^2f(\bbx_t)(\bbx_{t+1}-\bbx_t)-\nabla f(\bbx_{t+1})
 + \left(\tbH_t(K)-\bbH_t\right)(\bbx_{t+1}-\bbx_t).
\end{align}
\end{lemma}

\begin{myproof}
See Appendix \ref{app_lemma_esom}.
\end{myproof}

The results in Theorem \ref{lemma_esom} show the relationships between the primal $\bbx$ and dual $\bbv$ iterates generated by ESOM and the optimal arguments $\bbx^*$ and $\bbv^*$. The first result in \eqref{opt_res_ESOM1} is identical to the convergence property of PMM in \eqref{opt_res_pmm1}, while the second result in \eqref{opt_res_ESOM2} differs from \eqref{opt_res_pmm2} in having the extra summand $\bbe_t$. The vector $\bbe_t$ can be interpreted as the error of second order approximation for ESOM at step $t$. To be more precise, the optimality condition of the primal update of PMM is given by $\nabla f(\bbx_{t+1})+(\bbI-\bbZ)^{1/2}\bbv_t+\alpha(\bbI-\bbZ)\bbx_{t+1}+\eps (\bbx_{t+1}-\bbx_t)=\bb0$ as shown in \eqref{opt_res_pmm2}. Notice that the second order approximation of this condition is equivalent to $\nabla f(\bbx_{t})+\nabla^2 f(\bbx_{t})(\bbx_{t+1}-\bbx_t)+(\bbI-\bbZ)^{1/2}\bbv_t+\alpha(\bbI-\bbZ)\bbx_{t+1}+\eps (\bbx_{t+1}-\bbx_t)=\bb0$. However, the exact Hessian inverse $\bbH_t^{-1}=(\nabla^2 f(\bbx_{t})+\eps\bbI+\alpha(\bbI-\tbZ))^{-1}$ cannot be computed in a distributed manner to solve the optimality condition. Thus, it is approximated by the approximate Hessian inverse matrix $\tbH_t^{-1}(K)$ as introduced in \eqref{Hessian_approx}. This shows that the approximate optimality condition in ESOM is $\nabla f(\bbx_t)+ (\bbI-\bbZ)^{1/2}\bbv_t+\alpha(\bbI-\tbZ)\bbx_t+\tbH_t(\bbx_{t+1}-\bbx_{t})=\bb0$. Hence, the difference between the optimality conditions of PMM and ESOM is $\bbe_t=\nabla f(\bbx_t)-\nabla f(\bbx_{t+1})+\alpha(\bbI-\tbZ)(\bbx_t-\bbx_{t+1})+\tbH_t(\bbx_{t+1}-\bbx_{t})-\eps (\bbx_{t+1}-\bbx_t)$. By adding and subtracting the term $\bbH_t(\bbx_{t+1}-\bbx_t)$, the definition of the error vector $\bbe_t$ in \eqref{esom_error_vec} follows. 

The observation that the vector $\bbe_t$ characterizes the error of second order approximation in ESOM, motivates analyzing an upper bound for the error vector norm $\|\bbe_t\|$. To prove that the norm $\|\bbe_t\|$ is bounded above we assume the following condition is satisfied. 

%
\begin{assumption}\label{lip_hessian} 
The global objective function Hessian $\nabla^2 f(\bbx)$ is Lipschitz continuous with constant $L$, i.e.,
\begin{equation}\label{lip_ass}
\|\nabla^2 f(\bbx)-\nabla^2 f(\tbx)\|\leq L \|\bbx-\tbx\|.
\end{equation}
\end{assumption}

The conditions imposed by Assumption \ref{lip_hessian} is customary in the analysis of second order methods; see, e.g., \cite{mokhtari2015dqm}. In the following lemma we use the assumption in \eqref{lip_ass} to prove an upper bound for the error norm $\|\bbe_t\|$ in terms of the distance $\|\bbx_{t+1}-\bbx_t\|$.\\

\begin{lemma}\label{esom_approx_error}
Consider ESOM as introduced in \eqref{ESOM_primal_update}-\eqref{ESOM_dual_update} and recall the definition of the error vector $\bbe_t$ in \eqref{esom_error_vec}. Further, define $c>0$ as a lower bound for the local weights $\bbw_{ii}$. If Assumptions \ref{convexity_assumption}-\ref{lip_hessian} hold, then the error vector norm $\|\bbe_t\|$ is bounded above by
\begin{equation}\label{upp_bound_error}
\|\bbe_t\| \leq \Gamma_t \| \bbx_{t+1}-\bbx_t\|,
\end{equation}
where $\Gamma_t$ is defined as
\begin{equation}\label{gamma_def}
\Gamma_t:= \min\left\{2M, \frac{L}{2}\| \bbx_{t+1}-\bbx_t\|\right\}+\left({M+\eps+2\alpha(1-c)}\right)\rho^{K+1},
\end{equation}
and $\rho:={2\alpha(1-c)}/({2\alpha(1-c)+m+\eps})$.
\end{lemma}

\begin{myproof}
See Appendix \ref{app_lemma_esom_approx_error}.
\end{myproof}

First, note that the lower bound $c>0$ on the local weights $w_{ii}$ is implied from the fact that all the local weights are positive. In particular, we can define the lower bound $c$ as $c:=\min_{i}{w_{ii}}$. The result in \eqref{upp_bound_error} shows that the error of second order approximation in ESOM vanishes as the sequence of iterates $\bbx_t$ approaches the optimal argument $\bbx^*$. We will show in Theorem \ref{thm:esom_linear_convg} that $\|\bbx_{t}-\bbx^*\|$ converges to zero which implies that the limit of the sequence $ \|\bbx_{t+1}-\bbx_t\|$ is zero. 

To understand the definition of $\Gamma_t$ in \eqref{gamma_def}, we have to decompose the error vector $\bbe_t$ in \eqref{esom_error_vec} into two parts. The first part is $ \nabla f(\bbx_t)+\nabla^2f(\bbx_t)(\bbx_{t+1}-\bbx_t)-\nabla f(\bbx_{t+1})$ which comes from the fact that ESOM minimizes a second order approximation of the proximal augmented Lagrangian instead of the exact proximal augmented Lagrangian. This term can be bounded by $\min\{2M, ({L}/{2})\| \bbx_{t+1}-\bbx_t\|\}\|\bbx_{t+1}-\bbx_t\|$ as shown in Lemma \ref{esom_approx_error}. The second part of the error vector $\bbe_t$ is $(\tbH_t(K)-\bbH_t)(\bbx_{t+1}-\bbx_t)$ which shows the error of Hessian inverse approximation. Notice that computation of the exact Hessian inverse $\bbH_t^{-1}$ is not possible and ESOM approximates the exact Hessian by the approximation $\tbH_t^{-1}(K)$. According to the results in \cite{NN-part1}, the difference $\|\tbH_t(K)-\bbH_t\|$ can upper bounded by $(M+\eps+{2(1-c)}/{\alpha})\rho^{K+1}$ which justifies the second term of the expression for $\Gamma_t$ in \eqref{gamma_def}. In the following theorem, we use the result in Lemma \ref{esom_approx_error} to show that the sequence of Lyapunov functions $\|\bbu_t-\bbu^*\|_\bbG^2$ generated by ESOM converges to zero linearly.

\begin{thm}\label{thm:esom_linear_convg}
Consider ESOM as introduced in \eqref{ESOM_primal_update}-\eqref{ESOM_dual_update}. Consider $\beta>1$ and $\phi>1$ as arbitrary constants that are strictly larger than $1$, and $\zeta$ as a positive constant that is chosen from the interval $\zeta \in ( (m+M)/2mM, \eps/\Gamma_t^2 )$. Further, recall the definitions of the vector $\bbu$ and matrix $\bbG$ in \eqref{def:u_G} and consider $\tlmin(\bbI-\bbZ)$ as the smallest non-zero eigenvalue of the matrix $\bbI-\bbZ$. If Assumptions \ref{convexity_assumption}-\ref{lip_hessian} hold, then the sequence of Lyapunov functions $\|\bbu_{t}-\bbu^*\|_\bbG^2$ generated by ESOM satisfies 
\begin{equation}\label{ESOM_lin_convg}
\|\bbu_{t+1}-\bbu^*\|_\bbG^2 \ \leq\  \frac{1}{1+\delta_t'} \ \|\bbu_{t}-\bbu^*\|_\bbG^2.
\end{equation}
where the sequence $\delta'_t$ is given by
\begin{align}\label{esom_delta}
&\delta_t'=\min \Bigg\{ 	\frac{2\alpha \tlmin(\bbI-\bbZ)}{\phi\beta(m+M)},
	\left[\frac{2m M}{\eps(m+M)}\!-\!\frac{1}{\zeta\eps}\right],
	\frac{(\beta-1)\alpha\tlmin(\bbI-\bbZ)}{\beta\eps} \left[1-\frac{ \zeta\Gamma_t^2}{\eps}\right]  \! \left[1+\frac{\phi\Gamma_t^2(\beta-1) }{(\phi-1)\eps^2}\right]^{-1}\!\Bigg\}.
\end{align}
\end{thm}

\begin{myproof}
See Appendix \ref{app_thm_esom_linear_convg}.
\end{myproof}

The result in Theorem \ref{thm:esom_linear_convg} shows linear convergence of the sequence $\|\bbu_{t}-\bbu^*\|_\bbG^2$ generated by ESOM where the factor of linear convergence is $1/(1+\delta')$. Notice that the positive constant $\xi$ is chosen from the interval $ ( (m+M)/2mM, \eps/\Gamma_t^2 )$. This interval is non-empty if and only if the proximal parameter $\eps$ satisfies the condition $\eps>\Gamma_t^2(m+M)/2mM$. It follows from the result in Theorem \ref{thm:esom_linear_convg} that the sequence of primal variables $\bbx_t$ converges to the optimal argument $\bbx^*$ defined in \eqref{constrained_opt_problem}.\\

\begin{corollary}\label{esom_approx_error2}
Under the assumptions in Theorem \ref{thm:esom_linear_convg}, the sequence of squared errors $\|\bbx_t-\bbx^*\|^2$ generated by ESOM converges to zero at a linear rate, i.e.,
\begin{equation}\label{ESOM_lin_convg_22}
\|\bbx_{t}-\bbx^*\|^2  \leq  \left(\frac{1}{1+\min_{t}\{\delta_t'\}}\right)^t  \frac{\|\bbu_{0}-\bbu^*\|_\bbG^2}{\alpha \eps}.
\end{equation}
\end{corollary}
\begin{myproof}
According to the definition of the sequence $\bbu_t$ and matrix $\bbG$, we can write $\|\bbu_t-\bbu^*\|_\bbG^2=\alpha \eps \|\bbx_t-\bbx^*\|^2+ \|\bbv_t-\bbv^*\|^2$ which implies that $\|\bbx_t-\bbx^*\|^2\leq (1/\alpha \eps)\|\bbu_t-\bbu^*\|_\bbG^2$. Considering this result and linear convergence of the sequence $\|\bbu_t-\bbu^*\|_\bbG^2$ in \eqref{ESOM_lin_convg}, the claim in \eqref{ESOM_lin_convg_22} follows. 
\end{myproof}

\subsection{Convergence rates comparison}\label{sec:rate_comparison}

The expression for $\delta'_t$ in \eqref{esom_delta} verifies the intuition that the convergence rate of ESOM is slower than PMM. This is true, since the upper bounds for $\delta$ in PMM are larger than their equivalent upper bounds for $\delta'_t$ in ESOM. We obtain that $\delta_t'$ is smaller than $\delta$ which implies that the linear convergence factor $1/(1+\delta)$ of PMM is smaller than $1/(1+\delta'_t)$ for ESOM. Therefore, for all steps $t$, the linear convergence of PMM is faster than ESOM. Although, linear convergence factor of ESOM $1/(1+\delta'_t)$ is larger than $1/(1+\delta)$ for PMM, as time passes the gap between these two constants becomes smaller. In particular, notice that after a number of iterations $({L}/{2})\| \bbx_{t+1}-\bbx_t\|$ becomes smaller than $2M$ and $\Gamma_t$ can be simplified as 
\begin{align}\label{error_norm_bound2}
\Gamma_t \leq \frac{L}{2}\|\bbx_{t+1}-\bbx_t\|+ \left({2\alpha(1-c)}+M+\eps\right)\rho^{K+1}.
\end{align}
The term $(L/2)\| \bbx_{t+1}-\bbx_t\|$ eventually approaches zero, while the second term $({2(1-c)}/{\alpha}+M+\eps)\rho^{K+1}$ is constant. Although, the second term is not approaching zero, by proper choice of $\rho$ and $K$, this term can become arbitrary close to zero. Notice that when $\Gamma_t$ approaches zero, if we set $\zeta=1/\Gamma_t$ the upper bounds in \eqref{esom_delta} for $\delta_t'$ approach the upper bounds for $\delta$ of PMM in \eqref{pmm_delta}. 

Therefore, as time passes $\Gamma_t$ becomes smaller and the factor of linear convergence for ESOM $1/(1+\delta'_t)$ becomes closer to the linear convergence factor of PMM $1/(1+\delta)$.

\section{Numerical Experiments}\label{sec:numerical}

In this section we compare the performances of ESOM, EXTRA, Decentralized (D)ADMM, and Network Newton (NN).  First we consider a linear least squares problem and then we use the mentioned methods to solve a logistic regression problem.

\subsection{Decentralized linear least squares}\label{sec:DLS}
Consider a decentralized linear least squares problem where each agent
$i\in\{1,\cdots,n\}$ holds its private measurement equation,
$\bby_{i}=\bbM_{i} \tbx+\bbnu_{i}$, where $\bby_{i}\in\reals^{m_i}$ and $\bbM_{i}\in\reals^{m_i\times p}$ are measured data, $\tbx\in\reals^p$ is the unknown variable, and $\bbnu_{i}\in\mathbb{R}^{m_i}$ is some unknown noise. The decentralized linear least squares estimates $\tbx$ by solving the optimization problem 
\begin{equation}
\tbx^*=\argmin \limits_\tbx\ \sum\limits_{i=1}^n \|\bbM_{i} \tbx-\bby_{i}\|_2^2.
\end{equation}

The network in this experiment is randomly generated with
connectivity ratio $r={3}/{n}$, where $r$ is defined as the number of
edges divided by the number of all possible
ones, ${n(n-1)}/{2}$. We set $n=20$, $p=5$, and $m_i=5\ \forall\ i=1,\dots,n$. The vectors $\bby_{i}$ and matrices $\bbM_{i}$ as well as the noise vectors $\bbnu_{(i)}$, $\forall\ i$ are generated following the standard normal distribution. We precondition the aggregated data matrices $\bbM_{i}$ so that the condition number of the problem is $10$. The decision variables $\bbx_i$ are initialized as $\bbx_{i,0}=0$ for all nodes $i=1,\dots,n$ and the initial distance to the optimal is $\|\bbx_{i,0}-\bbx^*\|=100$. 

We use Metropolis constant edge weight matrix as the mixing matrix $\bbW$ in all experiments. We run PMM, EXTRA, and ESOM-$K$ with fixed hand-optimized stepsizes $\alpha$. The best choices of $\alpha$ for ESOM-0, ESOM-1, and ESOM-2 are $\alpha=0.03$, $\alpha=0.04$, and $\alpha=0.05$, respectively. The stepsize $\alpha=0.1$ leads to the best performance for EXTRA which is considered in the numerical experiments. Notice that for variations of NN-$K$, there is no optimal choice of stepsize -- smaller stepsize leads to more accurate but slow convergence, while large stepsize accelerates the convergence but to a less accurate neighborhood of the optimal solution. Therefore, for NN-0, NN-1, and NN-2 we set $\alpha=0.001$, $\alpha=0.008$, and $\alpha=0.02$, respectively. Although the PMM algorithm is not implementable in a decentralized fashion, we use its convergence path -- which is generated in a centralized manner -- as our benchmark. The choice of stepsize for PMM is $\alpha=2$.

\begin{figure}[t]
\begin{center}
\includegraphics[width=0.6\linewidth,height=0.4\linewidth]{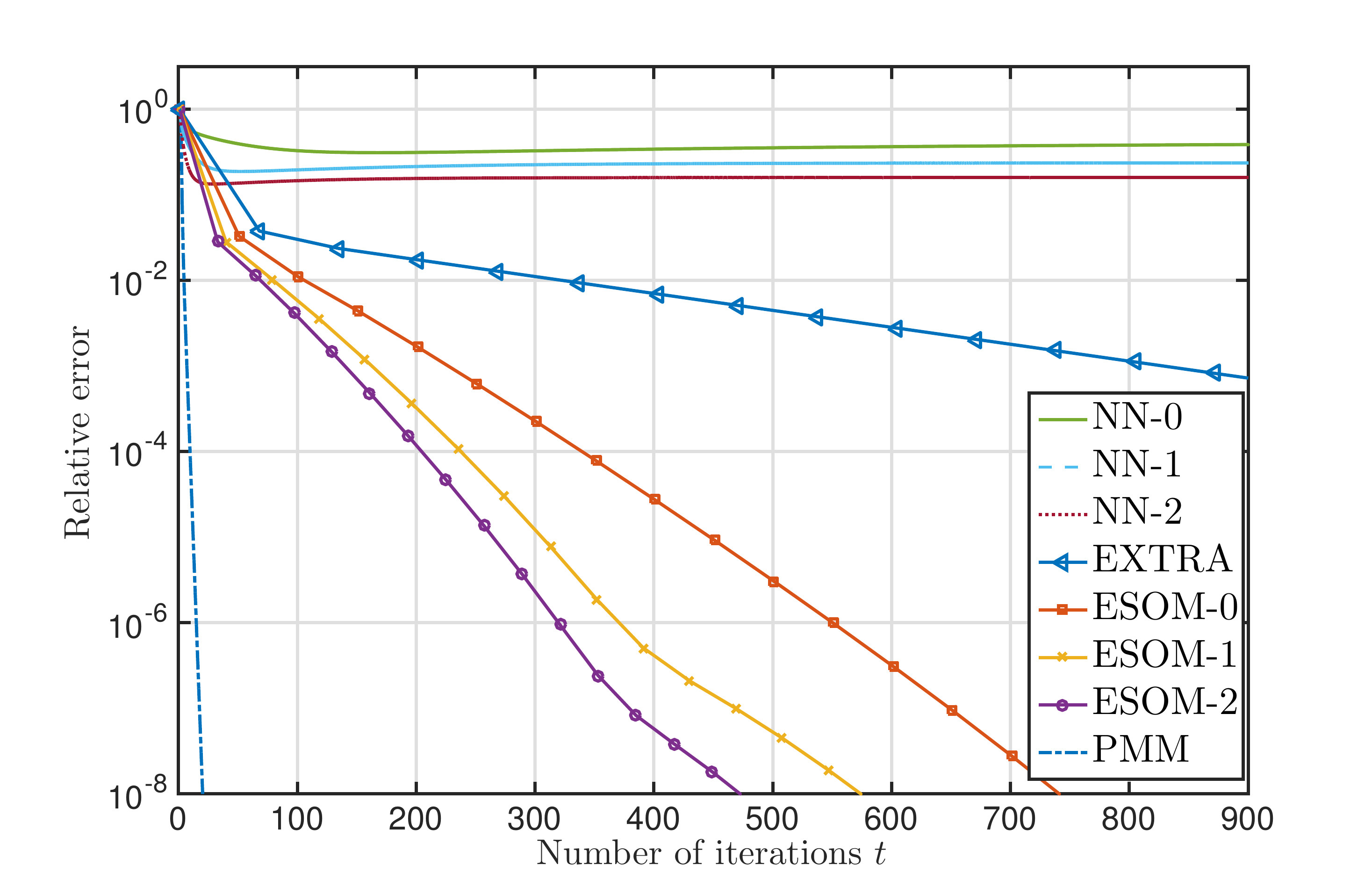}
\vspace{4mm}
\caption{Relative error ${\|\mathbf{x}_t-\mathbf{x}^*\|}/{\|\mathbf{x}_0-\mathbf{x}^*\|}$ of EXTRA, ESOM-$K$, NN-$K$, and PMM versus number of iterations for the least squares problem. Using a larger value of $K$ for ESOM-$K$ leads to faster convergence and makes the convergence path closer to the one for PMM. 
}\label{eps:num_LS_iter}
\end{center}
\end{figure}
\begin{figure}[t]
\begin{center}
\includegraphics[width=0.6\linewidth,height=0.4\linewidth]{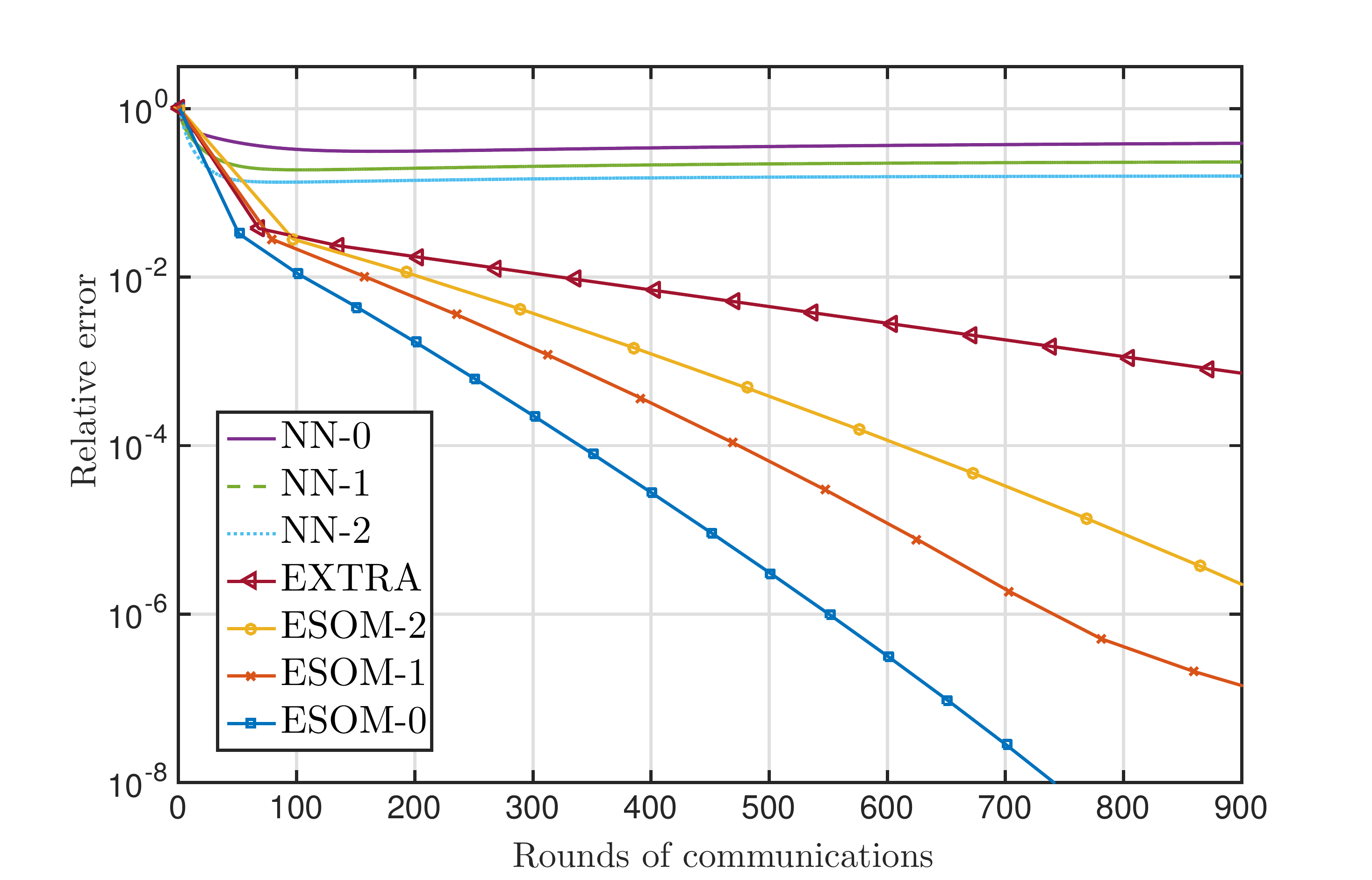}
\vspace{4mm}
\caption{Relative error ${\|\mathbf{x}_t-\mathbf{x}^*\|}/{\|\mathbf{x}_0-\mathbf{x}^*\|}$ of EXTRA, ESOM-$K$, NN-$K$, and PMM versus rounds of communications with neighboring nodes for the least squares problem. ESOM-$0$ is the most efficient algorithm in terms of communication cost among all the methods.}\label{eps:num_LS_comm}
\end{center}
\end{figure}

Fig. \ref{eps:num_LS_iter} illustrates the relative error ${\|\mathbf{x}_t-\mathbf{x}^*\|}/{\|\mathbf{x}_0-\mathbf{x}^*\|}$ versus the number of iterations. Notice that the vector $\bbx_t$ is the concatenation of the local vectors $\bbx_{i,t}$ and the optimal vector $\bbx^*$ is defined as $\bbx^*=[\tbx^*;\dots;\tbx^*]\in\reals^{np}$. Observe that all the variations of NN-$K$ fail to converge to the optimal argument and they converge linearly to a neighborhood of the optimal solution $\bbx^*$. Among the decentralized algorithms with exact linear convergence rate, EXTRA has the worst performance and all the variations of ESOM-$K$ outperform EXTRA. Recall that the problem condition number is $10$ in our experiment and the difference between EXTRA and ESOM-$K$ is more significant for problems with larger condition numbers. Further, choosing a larger value of $K$ for ESOM-$K$ leads to faster convergence and as we increase $K$ the convergence path of ESOM-$K$ approaches the convergence path of PMM.

EXTRA requires one round of communications per iteration, while NN-$K$ and ESOM-$K$ require $K+1$ rounds of local communications per iteration. Thus, convergence paths of these methods in terms of rounds of communications might be different from the ones in Fig. \ref{eps:num_LS_iter}. The convergence paths of NN, ESOM, EXTRA in terms of rounds of local communications are shown in Fig. \ref{eps:num_LS_comm}. In this plot we ignore PMM, since it requires infinite rounds of communications per iteration. The main difference between Figs. \ref{eps:num_LS_iter} and \ref{eps:num_LS_comm} is in the performances of ESOM-$0$, ESOM-$1$, and ESOM-$2$. All of the variations of ESOM outperform EXTRA in terms of rounds of communications, while the best performance belongs to ESOM-$0$. This observation shows that increasing the approximation level $K$ does not necessary improve the performance of ESOM-$K$ in terms of communication cost.

\subsection{Decentralized logistic regression}\label{sec:DLR}

We consider the application of ESOM for solving a logistic regression problem in a form
\begin{equation}
\tbx^*\! :=\argmin \limits_{\tbx\in \reals^p}  \frac{\lambda}{2}\|\tbx\|^2\!+\! \sum\limits_{i=1}^n\sum\limits_{j=1}^{m_i}
\ln\left(1+\exp\left(-(\bbs_{ij}^T\tbx)y_{ij}\right)\right),
\end{equation}
where every agent $i$ has access to $m_i$ training samples $\left(\bbs_{ij},y_{ij}\right)\in\R^p\times\{-1,+1\},\ j=1,\cdots,m_i$, including explanatory/feature variables $\bbs_{ij}$ and binary outputs/outcomes $y_{ij}$. The regularization term $(\lambda/2)\|\tbx\|^2$ is added to avoid overfitting where $\lambda$ is a positive constant. Hence, in the decentralized setting the local objective function $f_i$ of node $i$ is given by 
\begin{equation}
f_i(\tbx)=\frac{\lambda}{2n}\|\tbx\|^2\!+\! \sum\limits_{j=1}^{m_i}
\ln\left(1+\exp\left(-(\bbs_{ij}^T\tbx)y_{ij}\right)\right).
\end{equation}

The settings are as follows. The connected network is randomly generated with $n=20$ agents and connectivity ratio $r={3}/{n}$. Each agent holds $3$ samples, i.e., $m_i=3, \forall\ i$. The dimension of sample vectors $\bbs_{ij}$ is $p=3$. The samples are randomly generated, and the optimal logistic classifier $\tbx^*$ is pre-computed through centralized adaptive gradient method. We use Metropolis constant edge weight matrix as the mixing matrix $\bbW$ in ESOM-$K$. The stepsize $\alpha$ for ESOM-$0$, ESOM-$1$, ESOM-$2$, EXTRA, and DADMM are hand-optimized and the best of each is used for the comparison. 

\begin{figure}[t]
\begin{center}
\includegraphics[width=0.6\linewidth,height=0.4\linewidth]{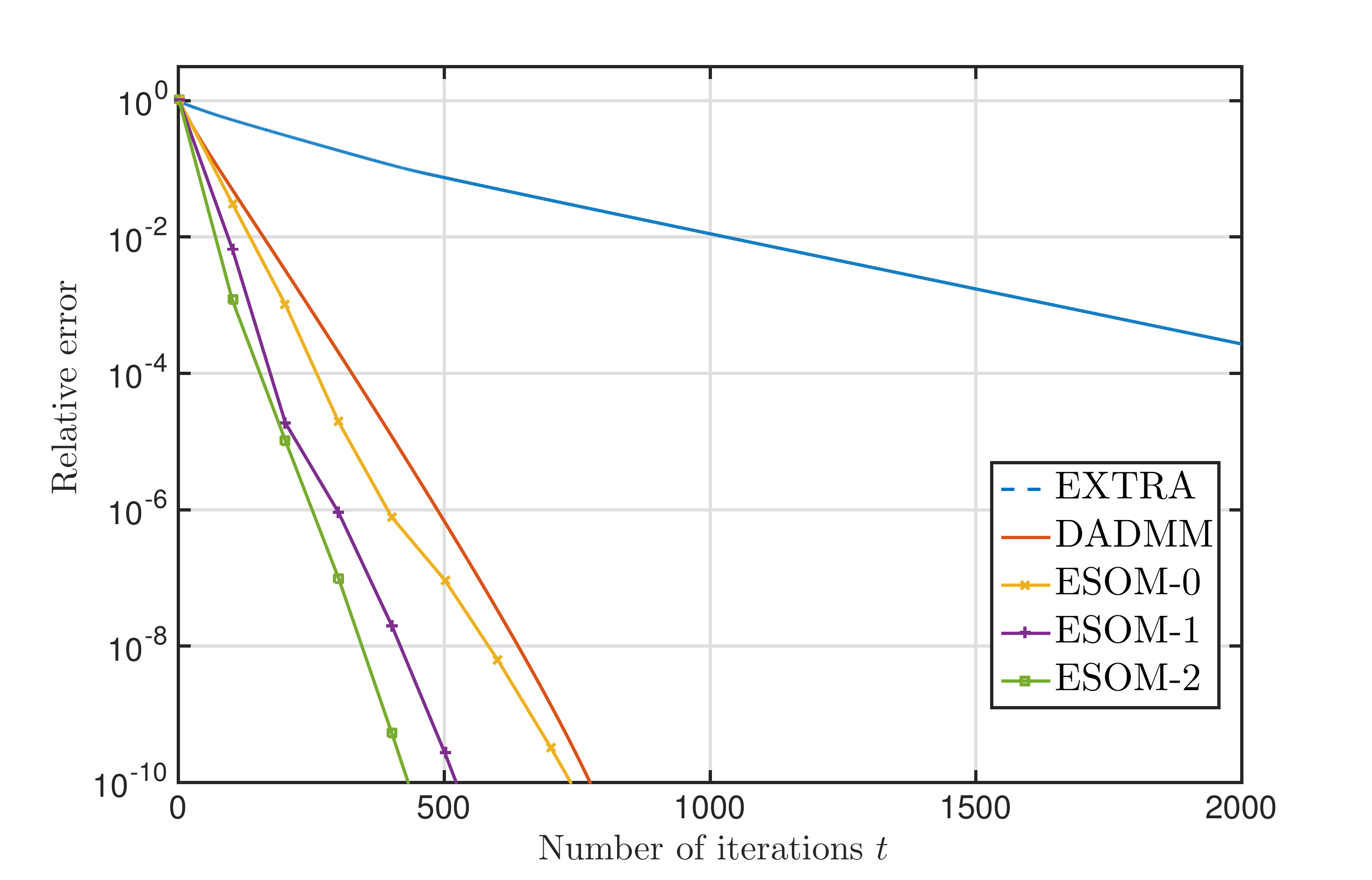}
\vspace{4mm}
\caption{Relative error ${\|\mathbf{x}_t-\mathbf{x}^*\|}/{\|\mathbf{x}_0-\mathbf{x}^*\|}$ of EXTRA, ESOM-$K$, and DADMM versus number of iterations for the logistic regression problem. EXTRA is significantly slower than the ESOM methods. The proposed methods (ESOM-$K$) outperform DADMM.}
\label{eps:num_Logistic_iter}
\end{center}
\end{figure}

\begin{figure}[t]
\begin{center}
\includegraphics[width=0.6\linewidth,height=0.4\linewidth]{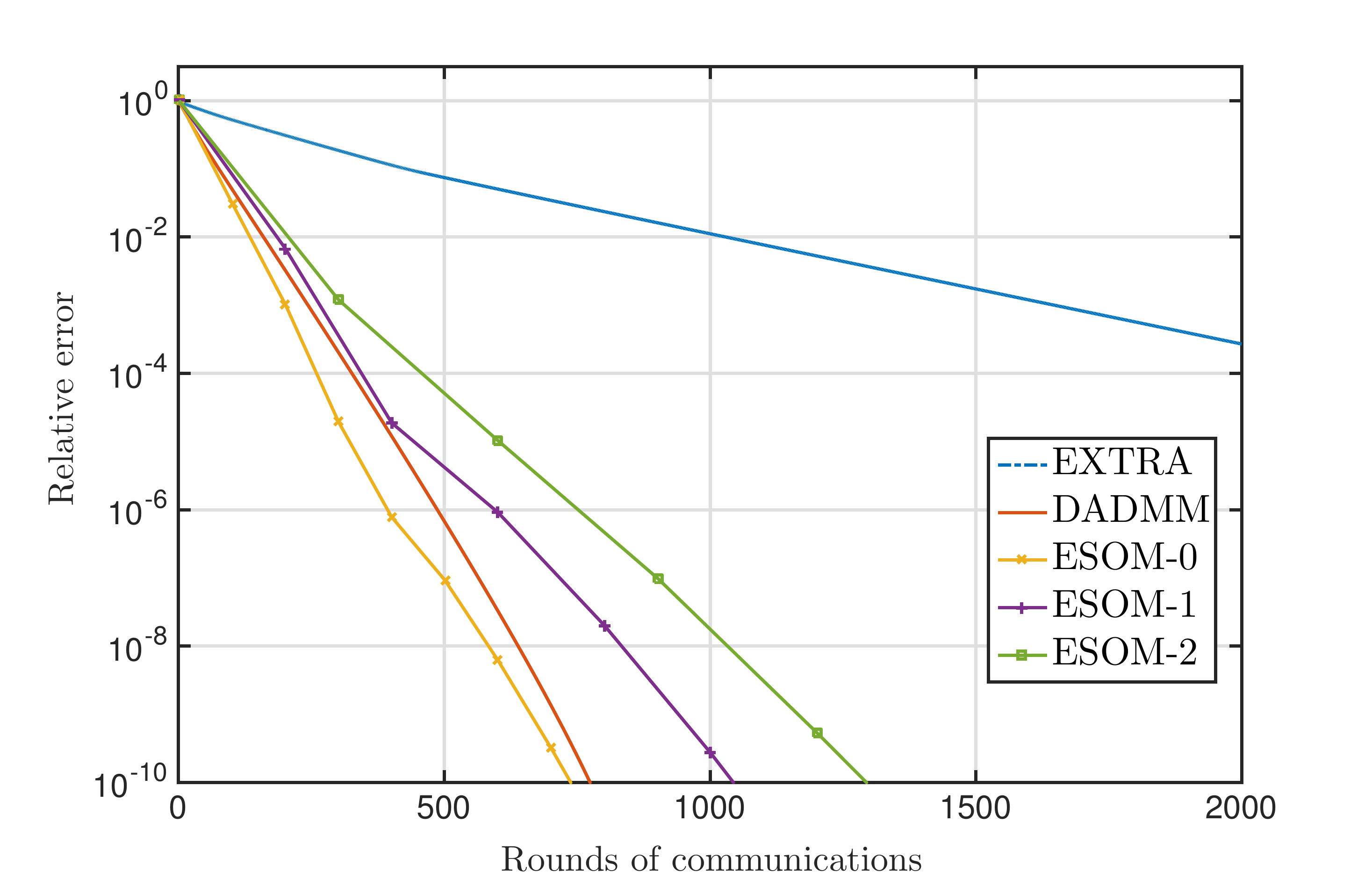}
\vspace{4mm}
\caption{Relative error ${\|\mathbf{x}_t-\mathbf{x}^*\|}/{\|\mathbf{x}_0-\mathbf{x}^*\|}$ of EXTRA, ESOM-$K$, and DADMM versus rounds of communications for the logistic regression problem. ESOM-$0$ has the best performance in terms of rounds of communications and it outperforms DADMM.}
\label{eps:num_Logistic_comm}
\end{center}
\end{figure}

Fig. \ref{eps:num_Logistic_iter} and Fig \ref{eps:num_Logistic_comm} showcase the convergence paths of ESOM-$0$, ESOM-$1$, ESOM-$2$, EXTRA, and DADMM versus number of iterations and rounds of communications, respectively. The results match the observations for the least squares problem in Fig. \ref{eps:num_LS_iter} and Fig. \ref{eps:num_LS_comm}. Different versions of ESOM-$K$ converge faster than EXTRA both in terms of communication cost and number of iterations. Moreover, ESOM-$2$ converges faster than ESOM-$1$ and ESOM-$0$ in terms of number of iterations, while ESOM-$0$ has the best performance in terms of communication cost for achieving a target accuracy. Comparing the convergence paths of ESOM-$0$, ESOM-$1$, and ESOM-$2$ with DADMM shows that number of iterations required for the convergence of DADMM is larger than the required iterations for ESOM-$0$, ESOM-$1$, and ESOM-$2$. In terms of communication cost, DADMM has a better performance relative to ESOM-$1$ and ESOM-$2$, while ESOM-$0$ is the  most efficient algorithm.


\section{Conclusions}\label{sec:conclusions}
We studied the consensus optimization problem where the components of a global objective function are available at different nodes of a network. We proposed an Exact Second-Order Method (ESOM) that converges to the optimal argument of the global objective function at a linear rate. We developed the update of ESOM by substituting the primal update of Proximal Method of Multipliers (PMM) with its second order approximation. Moreover, we approximated the Hessian inverse of the proximal augmented Lagrangian by truncating its Taylor's series. This approximation leads to a class of algorithms ESOM-$K$ where $K+1$ indicates the number of Taylor's series terms that are used for Hessian inverse approximation. Convergence analysis of ESOM-$K$ shows that the sequence of iterates converges to the optimal argument linearly irrespective to the choice of $K$. We showed that the linear convergence factor of ESOM-$K$ is a function of time and the choice of $K$. The linear convergence factor of ESOM approaches the linear convergence factor of PMM as time passes. Moreover, larger choice of $K$ makes the factor of linear convergence for ESOM closer to the one for PMM. Numerical results verify the theoretical linear convergence and the relation between the linear convergence factor of ESOM-$K$ and PMM. Further, we observed that larger choice of $K$ for ESOM-$K$ leads to faster convergence in terms of number of iterations, while the most efficient version of ESOM-$K$ in terms of communication cost is ESOM-$0$.

\begin{appendices}

%

\section{Proof of Lemma \ref{lemma_pmm}}\label{app_lemma_pmm}

Consider the updates of PMM in \eqref{PMM_primal_update} and \eqref{PMM_dual_update}. According to \eqref{constrained_opt_problem}, the optimal argument $\bbx^*$ satisfies the condition $(\bbI-\bbZ)^{1/2}\bbx^*=\bb0$. This observation in conjunction with the dual variable update in \eqref{PMM_dual_update} yields the claim in \eqref{opt_res_pmm1}.

To prove the claim in \eqref{opt_res_pmm2}, note that the optimality condition of \eqref{PMM_primal_update} implies $\nabla_\bbx \ccalL(\bbx_{t+1},\bbv_t)+\eps (\bbx_{t+1}-\bbx_t)=\bb0$. Based on the definition of the Lagrangian $\ccalL(\bbx,\bbv)$ in \eqref{lagrangian}, the optimality condition for the primal update of PMM can be written as
\begin{equation}\label{MM_primal_new_form}
\nabla f(\bbx_{t+1})+(\bbI-\bbZ)^{1/2}\bbv_t+\alpha(\bbI-\bbZ)\bbx_{t+1}+\eps (\bbx_{t+1}-\bbx_t)=\bb0.
\end{equation}
Further, notice that one of the KKT conditions of the optimization problem in \eqref{constrained_opt_problem} is  
\begin{equation}\label{KKT_condition}
\nabla f(\bbx^*) + (\bbI-\bbZ)^{1/2}\bbv^*=\bb0.
\end{equation}
Moreover, the optimal solution $\bbx^*=[\tbx^*;\dots;\tbx^*]$ of \eqref{constrained_opt_problem} lies in null$\{\bbI-\bbZ\}$. Therefore, we obtain 
\begin{equation}\label{opt_sol_consensus}
\alpha(\bbI-\bbZ)\bbx^*=\bb0.
\end{equation}
Subtracting the equalities in \eqref{KKT_condition} and \eqref{opt_sol_consensus} from \eqref{MM_primal_new_form} yields
\begin{align}\label{important_result_100}
\nabla f(\bbx_{t+1})-\nabla f(\bbx^*) +(\bbI-\bbZ)^{1/2}(\bbv_t-\bbv^*)
+\alpha(\bbI-\bbZ)(\bbx_{t+1}-\bbx^*)+\eps (\bbx_{t+1}-\bbx_t)=\bb0.
\end{align}
Regrouping the terms in \eqref{opt_res_pmm1} implies that $\bbv_t$ is equivalent to 
\begin{equation}\label{important_result_200}
\bbv_{t} =\bbv_{t+1}-\alpha (\bbI-\bbZ)^{1/2}(\bbx_{t+1}-\bbx^*).
\end{equation}
Substituting $\bbv_t$ in \eqref{important_result_100} by the expression in the right hand side of \eqref{important_result_200} follows the claim in \eqref{opt_res_pmm2}.

%
\section{Proof of Theorem \ref{thm:pmm_linear_convg}}\label{app_theorem_pmm}

According to Assumption \ref{convexity_assumption}, the global objective function $f$ is strongly convex with constant $m$ and its gradients $\nabla f$ are Lipschitz continuous with constant $M$. Considering these assumptions, we obtain that the inner product $(\bbx_{t+1}-\bbx^*)^T(\df(\bbx_{t+1})-\df(\bbx^*))$ is lower bounded by
\begin{align}\label{proof_001}
\frac{m M}{m+M}\|\bbx_{t+1}-\bbx^*\|^2+\frac{1}{m+M}\|\df(\bbx_{t+1})-\df(\bbx^*)\|^2 \leq (\bbx_{t+1}-\bbx^*)^T(\df(\bbx_{t+1})-\df(\bbx^*)).
\end{align}
The result in \eqref{opt_res_pmm2} shows that the difference $\df(\bbx_{t+1})-\df(\bbx^*)$ is equal to $-(\bbI-\bbZ)^{1/2}(\bbv_{t+1}-\bbv^*)-\eps (\bbx_{t+1}-\bbx_t)$. Apply this substitution into \eqref{proof_001} and multiply both sides of the resulted inequality by $2$ to obtain 
\begin{align}\label{proof_002}
&
\frac{2m M}{m+M}\|\bbx_{t+1}-\bbx^*\|^2+\frac{2}{m+M}\|\df(\bbx_{t+1})-\df(\bbx^*)\|^2\nonumber\\
&\quad \leq -2(\bbx_{t+1}-\bbx^*)^T(\bbI-\bbZ)^{1/2}(\bbv_{t+1}-\bbv^*)
-2\eps(\bbx_{t+1}-\bbx^*)^T(\bbx_{t+1}-\bbx_t).
\end{align}
Based on the result in \eqref{opt_res_pmm1}, we can substitute $(\bbx_{t+1}-\bbx^*)^T(\bbI-\bbZ)^{1/2}$ by $(1/\alpha)(\bbv_{t+1}-\bbv_t)^T$. Thus, we can rewrite \eqref{proof_002} as
\begin{align}\label{proof_003}
\!\!&
\frac{2 \alpha m M}{m+M}\|\bbx_{t+1}-\bbx^*\|^2+\frac{2\alpha }{m+M}\|\df(\bbx_{t+1})-\df(\bbx^*)\|^2\\
& \leq\! -2(\bbv_{t+1}-\bbv_t)^T\!(\bbv_{t+1}-\bbv^*)
 -2\alpha \eps(\bbx_{t+1}-\bbx^*)^T\!(\bbx_{t+1}-\bbx_t).\nonumber
\end{align}
Notice that for any vectors $\bba$, $\bbb$, and $\bbc$ we can write
\begin{equation}\label{proof_004}
2(\bba-\bbb)^T(\bba-\bbc)=\|\bba-\bbb\|^2+\|\bba-\bbc\|^2-\|\bbb-\bbc\|^2.
\end{equation}
By setting $\bba=\bbv_{t+1}$, $\bbb=\bbv_t$, and $\bbc=\bbv^*$ we obtain that the inner product $2(\bbv_{t+1}-\bbv_t)^T(\bbv_{t+1}-\bbv^*)$ in \eqref{proof_003} can be written as $\|\bbv_{t+1}-\bbv_t\|^2+\|\bbv_{t+1}-\bbv^*\|^2-\|\bbv_t-\bbv^*\|^2$. Likewise, setting $\bba=\bbx_{t+1}$, $\bbb=\bbx_t$, and $\bbc=\bbx^*$ implies that the inner product $2(\bbx_{t+1}-\bbx_t)^T(\bbx_{t+1}-\bbx^*)$ in \eqref{proof_003} is equal to $\|\bbx_{t+1}-\bbx_t\|^2+\|\bbx_{t+1}-\bbx^*\|^2-\|\bbx_t-\bbx^*\|^2$. Applying these simplifications into \eqref{proof_003} yields 
\begin{align}\label{proof_005}
&
\frac{2\alpha m M}{m+M}\|\bbx_{t+1}-\bbx^*\|^2+\frac{2\alpha}{m+M}\|\df(\bbx_{t+1})-\df(\bbx^*)\|^2\nonumber\\
& \leq\|\bbv_t-\bbv^*\|^2 -\|\bbv_{t+1}-\bbv_t\|^2-\|\bbv_{t+1}-\bbv^*\|^2+\alpha\eps\|\bbx_t-\bbx^*\|^2
 -\alpha\eps\|\bbx_{t+1}-\bbx_t\|^2-\alpha\eps\|\bbx_{t+1}-\bbx^*\|^2.
\end{align}
Now using the definitions of the variable $\bbu$ and matrix $\bbG$ in \eqref{def:u_G} we can substitute $\|\bbv_t-\bbv^*\|^2-\|\bbv_{t+1}-\bbv^*\|^2+\alpha\eps\|\bbx_t-\bbx^*\|^2-\alpha\eps\|\bbx_{t+1}-\bbx^*\|^2$ by $\|\bbu_{t}-\bbu^*\|_\bbG^2-\|\bbu_{t+1}-\bbu^*\|_\bbG^2$. Moreover, the squared norm $\|\bbv_{t+1}-\bbv_t\|^2$ is equivalent to $\|\bbx_{t+1}-\bbx^*\|_{\alpha^2(\bbI-\bbZ)}^2$ based on the result in \eqref{opt_res_pmm1}. By applying these substitutions we can rewrite \eqref{proof_005} as
\begin{align}\label{proof_006}
&
\frac{2\alpha m M}{m+M}\|\bbx_{t+1}-\bbx^*\|^2+\frac{2\alpha}{m+M}\|\df(\bbx_{t+1})-\df(\bbx^*)\|^2\nonumber\\
&\quad \leq\|\bbu_t-\bbu^*\|_\bbG^2 -\|\bbu_{t+1}-\bbu^*\|_\bbG^2
-\alpha\eps\|\bbx_{t+1}-\bbx_t\|^2
-\|\bbx_{t+1}-\bbx^*\|^2_{\alpha^2(\bbI-\bbZ)}.
\end{align}
Regrouping the terms in \eqref{proof_006} leads to the following lower bound for the difference $\|\bbu_t-\bbu^*\|_\bbG^2 -\|\bbu_{t+1}-\bbu^*\|_\bbG^2$,
\begin{align}\label{proof_007}
&
\|\bbu_t-\bbu^*\|_\bbG^2 -\|\bbu_{t+1}-\bbu^*\|_\bbG^2
\nonumber\\
&\quad\geq 
\frac{2\alpha}{m+M}\|\df(\bbx_{t+1})-\df(\bbx^*)\|^2
+\alpha\eps\|\bbx_{t+1}-\bbx_t\|^2
 +\|\bbx_{t+1}-\bbx^*\|^2_{\frac{2\alpha m M}{m+M}\bbI+\alpha^2(\bbI-\bbZ)}.
\end{align}
Observe that the result in \eqref{proof_007} provides a lower bound for the decrement $\|\bbu_t-\bbu^*\|_\bbG^2 -\|\bbu_{t+1}-\bbu^*\|_\bbG^2$. To prove the claim in \eqref{proof_0}, we need to show that for a positive constant $\delta$ we have $\|\bbu_t-\bbu^*\|_\bbG^2 -\|\bbu_{t+1}-\bbu^*\|_\bbG^2\geq \delta \|\bbu_{t+1}-\bbu^*\|_\bbG^2$. Therefore, the inequality in \eqref{proof_0} is satisfied if we can show that the lower bound in \eqref{proof_007} is greater than $\delta \|\bbu_{t+1}-\bbu^*\|_\bbG^2$ or equivalently 
\begin{align}\label{proof_008}
&\delta\|\bbv_{t+1}-\bbv^*\|^2+\delta\alpha\eps\|\bbx_{t+1}-\bbx^*\|^2
\nonumber\\
&
\leq 
\frac{2\alpha}{m+M}\|\df(\bbx_{t+1})-\df(\bbx^*)\|^2
+\alpha\eps\|\bbx_{t+1}-\bbx_t\|^2
+\|\bbx_{t+1}-\bbx^*\|^2_{\frac{2\alpha m M}{m+M}\bbI+\alpha^2(\bbI-\bbZ)}.
\end{align}
To prove that the inequality in \eqref{proof_008} for some $\delta>0$, we first find an upper bound for the squared norm $\|\bbv_{t+1}-\bbv^*\|^2$ in terms of the summands in the right hand side of \eqref{proof_008}. To do so, consider the relation \eqref{opt_res_pmm2} along with the fact that $\bbv_{t+1}$ and $\bbv^*$ both lying in the column space of $(\bbI-\bbZ)^{1/2}$. It follows that $\|\bbv_{t+1}-\bbv^*\|^2$ is bounded above by
\begin{align}\label{proof_009}
\|\bbv_{t+1}-\bbv^*\|^2
\leq
\frac{\beta\eps^2}{(\beta-1)\tlmin(\bbI-\bbZ)}\|\bbx_{t+1}-\bbx_t\|^2
 + \frac{\beta}{\tlmin(\bbI-\bbZ)}\|\df(\bbx_{t+1})-\df(\bbx^*)\|^2.
\end{align}
where $\beta>1$ is a tunable free parameter and $\tlmin(\bbI-\bbZ)$ is the smallest non-zero eigenvalue of $\bbI-\bbZ$. Considering the result in \eqref{proof_009} to satisfy the inequality in \eqref{proof_008}, which is a sufficient condition for the claim in \eqref{proof_0}, it remains to show that 
\begin{align}\label{proof_010}
&\frac{2\alpha}{m+M}\|\df(\bbx_{t+1})-\df(\bbx^*)\|^2
+\alpha\eps\|\bbx_{t+1}-\bbx_t\|^2
+\|\bbx_{t+1}-\bbx^*\|^2_{\frac{2\alpha m M}{m+M}\bbI+\alpha^2(\bbI-\bbZ)}\nonumber\\
&\geq   \frac{\delta\beta\eps^2}{(\beta-1)\tlmin(\bbI-\bbZ)}\|\bbx_{t+1}-\bbx_t\|^2
+\delta\eps\alpha\|\bbx_{t+1}-\bbx^*\|^2 + \frac{\delta\beta}{\tlmin(\bbI-\bbZ)}\|\df(\bbx_{t+1})-\df(\bbx^*)\|^2.
\end{align}
To enable \eqref{proof_010} and consequently enabling \eqref{proof_008}, we only need to verify that there exists $\delta>0$ such that
\begin{align}\label{cond_on_delta_pmm}
\frac{2\alpha m M}{m+M}\bbI+\alpha^2(\bbI-\bbZ) \succcurlyeq
\delta\alpha\eps\bbI, 
\quad \frac{2\alpha}{m+M}\geq \frac{\delta\beta}{\tlmin(\bbI-\bbZ)},
\quad \alpha \eps\geq\frac{\delta\beta\eps^2}{(\beta-1)\tlmin(\bbI-\bbZ)}.
\end{align}
The conditions in \eqref{cond_on_delta_pmm} are satisfied if the constant $\delta$ is chosen as in \eqref{pmm_delta}. Therefore, for $\delta$ in \eqref{pmm_delta} the claim in \eqref{proof_008} holds, which implies the claim in \eqref{proof_0}.

%
\section{Proof of Lemma \ref{lemma_esom}}\label{app_lemma_esom}

Consider the primal update of ESOM in \eqref{ESOM_primal_update_3}. By regrouping the terms we obtain that
\begin{align}\label{ESOM_gen_proof_001}
\nabla f(\bbx_t)+ (\bbI-\bbZ)^{1/2}\bbv_t+{\alpha}(\bbI-\bbZ)\bbx_t
 +\tbH_t(\bbx_{t+1}-\bbx_{t})=\bb0,
\end{align}
where $\tbH_t$ is the inverse of the Hessian inverse approximation $\tbH_t^{-1}(K)$. Recall the definition of the exact Hessian $\bbH_t$ in \eqref{exact_Hessian}. Adding and subtracting the term $\bbH_t(\bbx_{t+1}-\bbx_t)$ to the expression in \eqref{ESOM_gen_proof_001} yields
\begin{align}\label{ESOM_gen_proof_002}
\nabla f(\bbx_t)+\nabla^2f(\bbx_t)(\bbx_{t+1}-\bbx_t)+  (\bbI-\bbZ)^{1/2}\bbv_t
+{\alpha}(\bbI-\tbZ)\bbx_{t+1}
+{\eps}(\bbx_{t+1}-\bbx_t)+(\tbH_t-\bbH_t)(\bbx_{t+1}-\bbx_t)=\bb0.
\end{align}
Now using the definition of the error vector $\bbe_t$ in \eqref{esom_error_vec} we can rewrite \eqref{ESOM_gen_proof_002} as
\begin{align}\label{ESOM_gen_proof_003}
\nabla f(\bbx_{t+1})+(\bbI-\bbZ)^{1/2}\bbv_t
+{\alpha}(\bbI-\tbZ)\bbx_{t+1}
+{\eps}(\bbx_{t+1}-\bbx_t)+\bbe_t=\bb0.
\end{align}
Notice that the result in \eqref{ESOM_gen_proof_003} is identical to the expression for PMM in \eqref{MM_primal_new_form} except for the error term $\bbe_t$. To prove the claim in \eqref{opt_res_ESOM2} from \eqref{ESOM_gen_proof_003}, it remains to follow the steps in \eqref{KKT_condition}-\eqref{important_result_200}.

%
\section{Proof of Lemma \ref{esom_approx_error}}\label{app_lemma_esom_approx_error}

To prove the result in \eqref{upp_bound_error}, we first use the result in Proposition 2 of \cite{mokhtari2015dqm}. It shows that when the eigenvalues of the Hessian $\nabla^2 f(\bbx)$ are bounded above by $M$ and the Hessian is Lipschitz continuous with constant $L$ we can write
\begin{align}\label{proof_error_100}
\|\nabla f(\bbx_t)+\nabla^2f(\bbx_t)(\bbx_{t+1}-\bbx_t)-\nabla f(\bbx_{t+1})\|
\leq \min\left\{2M, \frac{L}{2}\| \bbx_{t+1}-\bbx_t\|\right\}.
\end{align}
Considering the result in \eqref{proof_error_100}, it remains to find an upper bound for the second term of the error vector $\bbe_t$ which is $(\tbH_t(K)-\bbH_t)(\bbx_{t+1}-\bbx_t)$. To do so, we develop first an upper bound for the norm $\|\tbH_t(K)-\bbH_t\|$. Notice that by factoring the term $\tbH_t(K)$ and using the Cauchy-Schwarz inequality we obtain that 
\begin{equation}\label{proof_error_200}
\left\|\tbH_t(K)-\bbH_t\right\|\leq \left\|\tbH_t(K)\right\|  \left\|\bbI-\bbH_t\tbH_t^{-1}(K)\right\|.
\end{equation}
According to Lemma 3 in \cite{NN-part1}, we can simplify $\bbI-\bbH_t\tbH_t^{-1}(K)$ as $(\bbB\bbD_t^{-1})^{K+1}$. This simplification implies that 
\begin{equation}\label{proof_error_300}
\left\|\bbI-\bbH_t\tbH_t^{-1}(K)\right\|=\left\|\bbB\bbD_t^{-1}\right\|^{K+1}.
\end{equation}
Observe that the matrices $\bbB$ and $\bbD_t$ in this paper are different from the ones in \cite{NN-part1}, but the analyses of them are very similar. Similar to the proof of Proposition 2 in \cite{NN-part1} we define $\hbD:=2\alpha(\bbI-\bbZ_d)$. Notice that the matrix $\hbD$ is bock diagonal where its $i$th diagonal block is $2\alpha(1-w_{ii})\bbI_p$. Thus, $\hbD$ is positive definite and invertible. Hence, We are allowed to write the product $\bbB\bbD_t^{-1}$ as
\begin{equation}\label{new_decomposition}
\bbB\bbD_t^{-1}=
	\left(\bbB\hbD^{-1}	\right)
	\left(	\hbD\bbD_t^{-1}	\right).
\end{equation}
The next step is to find an upper bound for the eigenvalues of $\bbB\hbD^{-1}$ in \eqref{new_decomposition}. Based on the definitions of matrices $\bbB$ and $\hbD$, the product $\bbB\hbD^{-1}$ is given by
\begin{equation}\label{product_of_B_and_D_hat}
\bbB\hbD^{-1} = \left(\bbI-2\bbZ_{d}+\bbZ \right) (2(\bbI-\bbZ_{d}))^{-1}.
\end{equation}
According to the result in Proposition 2 of \cite{NN-part1}, the eigenvalues of the matrix $(\bbI-2\bbZ_{d}+\bbZ ) (2(\bbI-\bbZ_{d}))^{-1}$ are uniformly bounded by $0$ and $1$. Thus, we obtain that
 \begin{equation}\label{kosenanat}
 \|\bbB\hbD^{-1}\|\leq 1.
\end{equation}

According to the definitions of the matrices $\hbD$ and $\bbD_t$, the product $\hbD^{1/2}\bbD_t^{-1/2}$ is block diagonal and the $i$th diagonal block is given by
\begin{equation}\label{product_DDD}
\left[\hbD\bbD_{t}^{-1}\right]_{ii}= 
\left(\frac{  \nabla^2 f_{i}(\bbx_{i,t})+\eps \bbI}{2\alpha(1-w_{ii})} +\bbI\right)^{-1}.
\end{equation}
Based on Assumption \ref{convexity_assumption}, the eigenvalues of the local Hessians $\nabla^{2}f_{i}(\bbx_{i})$ are bounded by $m$ and $M$. Further, notice that the diagonal elements $w_{ii}$ of the weight matrix $\bbW$ are bounded below by $c$. Considering these bounds, we can show that the eigenvalues of the matrices $(1/ 2\alpha(1-w_{ii})) (\nabla^2 f_{i}(\bbx_{i,t})+\eps\bbI)+\bbI$ for all $i=1,\dots,n$ are bounded below by 
\begin{equation}\label{bounds_33}
\left[\frac{ m+\eps}{2\alpha(1-c)}+1\right]  \bbI  
\preceq  \frac{  \nabla^2 f_{i}(\bbx_{i,t})+\eps \bbI}{2\alpha(1-w_{ii})} +\bbI .
\end{equation}
By considering the bounds in \eqref{bounds_33}, the eigenvalues of each block of the matrix $\hbD\bbD_t^{-1}$, introduced in \eqref{product_DDD}, are bounded above as
\begin{align}\label{bounds_3}
\left(\frac{  \nabla^2 f_{i}(\bbx_{i,t})+\eps \bbI}{2\alpha(1-w_{ii})} +\bbI\right)^{-1} \preceq  \left[\frac{ m+\eps}{2\alpha(1-c)}+1\right]^{-1}  \bbI.
\end{align}
The upper bound in \eqref{bounds_3} for the eigenvalues of each diagonal block of the matrix $\hbD\bbD_t^{-1}$ implies that the matrix norm $\|\hbD\bbD_t^{-1}\|$ is bounded above by 
\begin{align}\label{bounds_390900}
\|\hbD\bbD_t^{-1}\| \leq \rho:= \frac{2\alpha(1-c)}{2\alpha(1-c)+m+\eps}.
\end{align}
Considering the upper bounds in \eqref{kosenanat} and \eqref{bounds_390900} and the relation in \eqref{new_decomposition} we obtain that 
\begin{align}\label{bounds_39090}
\|\bbB\bbD_t^{-1}\| \leq \rho.
\end{align}
Substituting the norm $\|\bbB\bbD_t^{-1}\|$ in \eqref{proof_error_300} by its upper bound $\rho$ implies $\|\bbI-\bbH_t\tbH_t^{-1}(K)\|\leq \rho^{K+1} $. This result in conjunction with the inequality in \eqref{proof_error_200} yields 
\begin{equation}\label{proof_error_500}
\left\|\tbH_t(K)-\bbH_t\right\|\leq \rho^{K+1} \left\|\tbH_t(K)\right\| .
\end{equation}

To bound the norm $\|\tbH_t(K)\| $, we first find a lower bound for the eigenvalues of the approximate Hessian inverse $\tbH_t^{-1}(K)$. Notice that according to the definition of the approximate Hessian inverse in \eqref{Hessian_approx}, we can write 
\begin{equation}\label{Hessian_approx222}
\tbH_t^{-1}(K):=\bbD_t^{-1}+\bbD_t^{-1}\
\sum_{u=1}^K(\bbD_t^{-1/2}\bbB\bbD_t^{-1/2})^u\
\bbD_t^{-1/2}.
\end{equation}
Notice that according to the result in Proposition 1 of \cite{NN-part1}, the matrix $ \left(\bbI-2\bbZ_{d}+\bbZ \right) $ is positive semidefinite which implies that $\bbB=\alpha\left(\bbI-2\bbZ_{d}+\bbZ \right) $ is also positive semidefinite. Thus, all the $K$ summands in \eqref{Hessian_approx222} are positive semidefinite and as a result we obtain that
\begin{equation}\label{lower_bounded_matrix}
\bbD_t^{-1} \preceq\  \tbH_t^{-1}(K) .
\end{equation}
The eigenvalues of $\bbI-\bbZ_d$ are bounded above by $1-c$, since all the local weights $w_{ii}$ are larger than $c$. This observation in conjunction with the strong convexity of the global objective function $f$ implies that the eigenvalues of $\bbD_t = \nabla^2f(\bbx_t)+\eps\bbI+2\alpha(\bbI-\bbZ_d)$ are bounded above by $M+\eps+2\alpha(1-c)$. Therefore,
\begin{equation}\label{lower_bounded_matrix2}
\frac{1}{M+\eps+2\alpha(1-c)}\ \bbI \ \preceq\ \bbD_t^{-1}.
\end{equation}
The results in \eqref{lower_bounded_matrix} and \eqref{lower_bounded_matrix2} imply that the eigenvalues of the approximate Hessian inverse $\tbH_t^{-1}(K)$ are greater than ${1}/({M+\eps+2\alpha(1-c)})$. Therefore, the eigenvalues of the positive definite matrix $\tbH_t(K)$ are smaller than ${M+\eps+2\alpha(1-c)}$ and we can write 
\begin{equation}\label{proof_error_700}
\left\|\tbH_t(K)\right\| \leq {M+\eps+2\alpha(1-c)}. 
\end{equation}
Considering the inequalities in \eqref{proof_error_500} and \eqref{proof_error_700} and using the Cauchy-Schwarz inequality we can show that the norm $\|(\tbH_t(K)-\bbH_t)(\bbx_{t+1}-\bbx_t)\|$ is bounded above by 
\begin{align}\label{proof_error_800}
\left\|(\tbH_t(K)-\bbH_t)(\bbx_{t+1}-\bbx_t)\right\| 
\leq \left({M+\eps+2\alpha(1-c)}\right)\rho^{K+1} \|\bbx_{t+1}-\bbx_t\|.\qquad
\end{align}
Observing the inequalities in \eqref{proof_error_100} and \eqref{proof_error_800} and using the triangle inequality the claim in \eqref{upp_bound_error} follows.

%
\section{Proof of Theorem \ref{thm:esom_linear_convg}}\label{app_thm_esom_linear_convg}

Notice that in proving the claim in \eqref{ESOM_lin_convg} we use some of the steps in the proof of Theorem \ref{thm:pmm_linear_convg} to avoid rewriting similar equations. First, note that according to the result in \eqref{opt_res_ESOM2}, the difference $\nabla f(\bbx_{t+1})-\nabla f(\bbx^*)$ for the ESOM method can be written as 
\begin{align}\label{proof_lin_ESOM_001}
\nabla f(\bbx_{t+1})-\nabla f(\bbx^*) =-(\bbI-\bbZ)^{1/2}(\bbv_{t+1}-\bbv^*)-\eps (\bbx_{t+1}-\bbx_t)-\bbe_t.
\end{align}
Now recall the the inequality in \eqref{proof_001} and substitute the gradients difference $\nabla f(\bbx_{t+1})-\nabla f(\bbx^*) $ in the inner product $(\bbx_{t+1}-\bbx^*)^T(\df(\bbx_{t+1})-\df(\bbx^*))$ by the expression in the right hand side of \eqref{proof_lin_ESOM_001}. Applying this substitution and multiplying both sides of the implied inequality by $2\alpha$ follows 
\begin{align}\label{proof_lin_ESOM_002}
&
\frac{2\alpha m M}{m+M}\|\bbx_{t+1}-\bbx^*\|^2+\frac{2\alpha}{m+M}\|\df(\bbx_{t+1})-\df(\bbx^*)\|^2\nonumber\\
& \leq -2\alpha(\bbx_{t+1}-\bbx^*)^T(\bbI-\bbZ)^{1/2}(\bbv_{t+1}-\bbv^*)
-2\alpha\eps(\bbx_{t+1}-\bbx^*)^T(\bbx_{t+1}-\bbx_t)
-2\alpha(\bbx_{t+1}-\bbx^*)^T\bbe_t.
\end{align}
By following the steps in \eqref{proof_002}-\eqref{proof_007}, the result in \eqref{proof_lin_ESOM_002} leads to a lower bound for the difference $\|\bbu_t-\bbu^*\|_\bbG^2 -\|\bbu_{t+1}-\bbu^*\|_\bbG^2$ as 
\begin{align}\label{proof_lin_ESOM_003}
&
\|\bbu_t-\bbu^*\|_\bbG^2 -\|\bbu_{t+1}-\bbu^*\|_\bbG^2
\nonumber\\
&\geq 
\frac{2\alpha}{m+M}\|\df(\bbx_{t+1})-\df(\bbx^*)\|^2
+\alpha\eps\|\bbx_{t+1}-\bbx_t\|^2
+\|\bbx_{t+1}-\bbx^*\|^2_{\frac{2\alpha m M}{m+M}\bbI+\alpha^2(\bbI-\bbZ)}+2\alpha(\bbx_{t+1}-\bbx^*)^T\bbe_t.
\end{align}
Notice that the inner product $2(\bbx_{t+1}-\bbx^*)^T\bbe_t$ is bounded below by $-(1/\zeta)\|\bbx_{t+1}-\bbx^*\|^2-\zeta \|\bbe_t\|^2$ for any positive constant $\zeta>0$. Therefore, the lower bound in \eqref{proof_lin_ESOM_003} can be updated as
\begin{align}\label{proof_lin_ESOM_004}
&\|\bbu_t-\bbu^*\|_\bbG^2-\|\bbu_{t+1}-\bbu^*\|_\bbG^2\nonumber\\
& \geq
\|\bbx_{t+1}-\bbx^*\|_{(\frac{2\alpha m M}{m+M}-\frac{\alpha}{\zeta})\bbI+\alpha^2(\bbI-\bbZ)}^2+\alpha\eps\|\bbx_{t+1}-\bbx_{t}\|^2
+\frac{2\alpha}{m+M}\|\df(\bbx_{t+1})-\df(\bbx^*)\|_\Fro^2-\alpha \zeta\|\bbe_t\|^2.
\end{align}
In order to establish \eqref{ESOM_lin_convg}, we need to show that the difference $\|\bbu_t-\bbu^*\|_\bbG^2-\|\bbu_{t+1}-\bbu^*\|_\bbG^2$ is bounded below by $\delta_t'\|\bbu_{t+1}-\bbu^*\|_\bbG^2$. To do so, we show that the lower bound for $\|\bbu_t-\bbu^*\|_\bbG^2-\|\bbu_{t+1}-\bbu^*\|_\bbG^2$ in \eqref{proof_lin_ESOM_004} is larger than $\delta_t'\|\bbu_{t+1}-\bbu^*\|_\bbG^2$, i.e., 
\begin{align}\label{proof_lin_ESOM_005}
&\delta_t'\|\bbv_{t+1}-\bbv^*\|^2+\delta_t'\alpha\eps\|\bbx_{t+1}-\bbx^*\|^2\nonumber\\
& \leq
\|\bbx_{t+1}-\bbx^*\|_{(\frac{2\alpha m M}{m+M}-\frac{\alpha}{\zeta})\bbI+\alpha^2(\bbI-\bbZ)}^2+\alpha\eps\|\bbx_{t+1}-\bbx_{t}\|^2
 +\frac{2\alpha}{m+M}\|\df(\bbx_{t+1})-\df(\bbx^*)\|_\Fro^2-\alpha \zeta\|\bbe_t\|^2.
\end{align}
We proceed to find an upper bound for the squared norm $\|\bbv_{t+1}-\bbv^*\|^2$ in terms of the summands in the right hand side of \eqref{proof_lin_ESOM_005}. Consider the relation \eqref{ESOM_gen_proof_003} as well as the fact that $\bbv_{t+1}$ and $\bbv^*$ both lie in the column space of $(\bbI-\bbZ)^{1/2}$. It follows that $\|\bbv_{t+1}-\bbv^*\|^2$ is bounded above by
\begin{align}\label{proof_lin_ESOM_006}
\|\bbv_{t+1}-\bbv^*\|^2
&\leq
\frac{\beta\eps^2}{(\beta-1)\tlmin(\bbI-\bbZ)}\|\bbx_{t+1}-\bbx_t\|^2
+ \frac{\phi\beta}{\tlmin(\bbI-\bbZ)}\|\df(\bbx_{t+1})-\df(\bbx^*)\|^2
\nonumber\\
&  \quad
+\frac{\beta\phi}{(\phi-1)\tlmin(\bbI-\bbZ)}\|\bbe_t\|^2 .
\end{align}
By substituting the upper bound in \eqref{proof_lin_ESOM_006} for the squared norm $\|\bbv_{t+1}-\bbv^*\|^2$ in \eqref{proof_lin_ESOM_005} we obtain a sufficient condition for the result in \eqref{proof_lin_ESOM_005} which is given by 
\begin{align}\label{proof_lin_ESOM_007}
&\delta_t'\alpha\eps\|\bbx_{t+1}-\bbx^*\|^2
+
\frac{\delta'\beta\eps^2}{(\beta-1)\tlmin(\bbI-\bbZ)}\|\bbx_{t+1}-\bbx_t\|^2
\nonumber\\
&+ \frac{\delta_t'\phi\beta}{\tlmin(\bbI-\bbZ)}\|\df(\bbx_{t+1})-\df(\bbx^*)\|^2
+\frac{\delta_t'\beta\phi\alpha^2\|\bbe_t\|^2 }{(\phi-1)\tlmin(\bbI-\bbZ)}\nonumber\\
&\quad  \leq
\|\bbx_{t+1}-\bbx^*\|_{(\frac{2\alpha m M}{m+M}-\frac{\alpha}{\zeta})\bbI+\alpha^2(\bbI-\bbZ)}^2+\alpha\eps\|\bbx_{t+1}-\bbx_{t}\|^2
 +\frac{2\alpha}{m+M}\|\df(\bbx_{t+1})-\df(\bbx^*)\|^2-\alpha \zeta\|\bbe_t\|^2.
\end{align}
Substitute the squared norm $\|\bbe_t\|^2$ terms in \eqref{proof_lin_ESOM_007} by the upper bound in \eqref{upp_bound_error}. It follows from this substitution and regrouping the terms that 
\begin{align}\label{proof_lin_ESOM_008}
&  
0 \leq
\|\bbx_{t+1}-\bbx^*\|_{(\frac{2\alpha m M}{m+M}-\frac{\alpha}{\zeta}-\delta_t'\alpha\eps)\bbI+\alpha^2(\bbI-\bbZ)}^2
+\left(\frac{2\alpha}{m+M}-\frac{\delta_t'\phi\beta}{\tlmin(\bbI-\bbZ)}\right)\|\df(\bbx_{t+1})-\df(\bbx^*)\|^2
\nonumber\\
&\quad
+\Bigg[\alpha\eps-\frac{\delta_t'\beta\eps^2}{(\beta-1)\tlmin(\bbI-\bbZ)}
-\frac{\delta_t'\beta\phi\Gamma^2  }{(\phi-1)\tlmin(\bbI-\bbZ)}-\alpha \zeta\Gamma^2  \Bigg]\|\bbx_{t+1}-\bbx_{t}\|^2.
\end{align}
Notice that if the inequality in \eqref{proof_lin_ESOM_008} is satisfied, then the result in \eqref{proof_lin_ESOM_007} holds which implies the result in \eqref{proof_lin_ESOM_005} and the linear convergence claim in \eqref{ESOM_lin_convg}. To satisfy the inequality in \eqref{proof_lin_ESOM_008} we need to make sure that the coefficients of the terms $\|\bbx_{t+1}-\bbx_{t}\|^2$, $\|\bbx_{t+1}-\bbx^*\|^2$, and $\|\df(\bbx_{t+1})-\df(\bbx^*)\|^2$ are non-negative. Therefore, the inequality in  \eqref{proof_lin_ESOM_008} holds if $\delta_t'$ satisfies 
\begin{align}\label{havij}
&\frac{2\alpha m M}{m+M}-\frac{\alpha}{\zeta}-\delta_t'\alpha\eps\geq 0,
\quad \frac{2\alpha}{m+M}\geq \frac{\delta_t'\phi\beta}{\tlmin(\bbI-\bbZ)}\\
&\alpha \eps\geq\frac{\delta_t'\beta\eps^2}{(\beta-1)\tlmin(\bbI-\bbZ)}+\frac{\delta_t'\beta\phi\Gamma^2  }{(\phi-1)\tlmin(\bbI-\bbZ)}+\alpha \zeta\Gamma^2.\nonumber
\end{align}
The conditions in \eqref{havij} are satisfied if $\delta_t'$ is chosen as in \eqref{esom_delta}. Thus, $\delta_t'$ in \eqref{esom_delta} satisfies the conditions in \eqref{havij} and the claim in \eqref{ESOM_lin_convg} holds. 
\end{appendices}

\bibliographystyle{IEEEtran}
  \bibliography{bmc_article}
\end{document}